\documentclass[draftcls,onecolumn]{IEEEtran}
\usepackage{cite}
\usepackage{amsmath,amssymb,amsfonts}
\usepackage{algorithmic}
\usepackage{graphicx}
\usepackage{textcomp}
\usepackage{float}
\usepackage[subrefformat=parens,labelformat=parens,caption=false]{subfig}
\usepackage[flushleft]{threeparttable}
\allowdisplaybreaks
\usepackage{xcolor}

\def\BibTeX{{\rm B\kern-.05em{\sc i\kern-.025em b}\kern-.08em
    T\kern-.1667em\lower.7ex\hbox{E}\kern-.125emX}}    
\begin{document}
\title{Multiphysics Simulation of Plasmonic Photoconductive Devices using Discontinuous Galerkin Methods}

\author{Liang Chen and Hakan Bagci
\thanks{The authors are with the Division of Computer, Electrical, and Mathematical Science and Engineering, King Abdullah University of Science and Technology (KAUST), Thuwal 23955-6900, Saudi Arabia (e-mails:\{liang.chen, hakan.bagci\}@kaust.edu.sa).}}
\maketitle

\begin{abstract}
Plasmonic nanostructures significantly improve the performance of photoconductive devices (PCDs) in generating terahertz radiation. However, they are geometrically intricate and result in complicated electromagnetic (EM) field and carrier interactions under a bias voltage and upon excitation by an optical EM wave. These lead to new challenges in simulations of plasmonic PCDs, which cannot be addressed by existing numerical frameworks. In this work, a multiphysics framework making use of discontinuous Galerkin (DG) methods is developed to address these challenges. The operation of the PCD is analyzed in stationary and transient states, which are described by coupled systems of the Poisson and stationary drift-diffusion (DD) equations and the time-dependent Maxwell and DD equations, respectively. Both systems are discretized using DG schemes. The nonlinearity of the stationary system is accounted for using the Gummel iterative method while the nonlinear coupling between the time-dependent Maxwell and DD equations is tackled during time integration. The DG-based discretization and the explicit time marching help in handling space and time characteristic scales that are associated with different physical processes and differ by several orders of magnitude. The accuracy and applicability of the resulting multiphysics framework are demonstrated via simulations of conventional and plasmonic PCDs.
\end{abstract}


\section{Introduction}
\IEEEPARstart{I}{n} recent years, electromagnetic (EM) systems operating at terahertz (THz) frequencies have found applications in various fields ranging from wireless communication to non-destructive testing~\cite{Federici2010review,Zaytsev2019review}. One of the fundamental challenges in developing these systems is implementing efficient THz sources and detectors. Plasmonic photoconductive devices (PCDs) have become one of the most promising candidates to address this challenge because of their dramatically increased optical-to-THz conversion efficiency~\cite{Lepeshov2017review, Burford2017review, Kang2018review, Yardimci2018review, Yachmenev2019review}. While the experimental research on fabrication of plasmonic PCDs has progressed remarkably in the recent years, there is still room for significant improvement in formulating and implementing numerical methods for rigorous and accurate simulation of this type of devices.

Challenges in simulation of PCDs stem from the fact that their operation relies on strongly coupled and simultaneously occurring multiple physical processes. A PCD incorporates a photoconductive semiconductor device that is biased under an external applied voltage. Upon excitation by an optical EM wave, this semiconductor device generates carriers and produces THz currents~\cite{Lepeshov2017review, Burford2017review}. Therefore, the operation of the PCD can be analyzed in two stages: (1) a stationary state corresponding to the non-equilibrium state of the semiconductor device under the bias voltage, and (2) a transient stage describing the dynamic processes after the optical EM wave impinges on the device. The stationary state is a result of the interaction between the static electric field and the carriers and is mathematically described by a coupled system of the Poisson equation and the stationary drift-diffusion (DD) equations. At the transient stage, several coupled processes happen simultaneously: EM field/wave interactions (propagation and scattering), carrier generation and recombination, and carrier drift and diffusion. This stage is mathematically described by a coupled system of the time-dependent Maxwell and DD equations. Because of the strong nonlinear coupling between these two sets of equations, 
accurate modeling of the transient stage calls for a multiphysics solver operating in time-domain.

Characteristic space and time scales of the various physical interactions listed above differ by several orders of magnitude. Space scales change from $\sim  10~\mathrm{nm}$ (Debye length) to $\sim  1~\mu$m (geometrical dimensions of the device). Furthermore, plasmonic nanostructures and localized plasmon modes generated on them make the space scale even more complicated. Time scales are determined by the carrier dynamics (including the advection and diffusion motions) and the EM wave interactions. These two scales typically differ by $\sim 3$ orders of magnitude. This multiscale nature of the physical interactions involved in the PCD operation requires highly flexible spatial and temporal discretization techniques. 

Several numerical methods have been developed for simulating PCDs. For conventional PCDs (PCDs without plasmonic nanostructures), the finite difference time-domain (FDTD) method has been extensively used~\cite{Kirawanich2008, Moreno2014, Moreno20141}. These schemes solve the time-dependent Maxwell and DD equations to account for carrier dynamics and the THz EM wave radiation of the attached THz antenna, respectively, while the carrier generation due to the optical EM wave excitation is taken into account using a closed-form analytical expression (the optical EM wave interactions are not simulated). It has been shown in~\cite{Moreno2014, Moreno20141} that the results obtained using this approach match experimental results well for conventional PCDs. For a plasmonic PCD, the nanostructures introduced on top of or inside the semiconductor region  produce highly localized and strong EM fields resulting in a significantly increased carrier generation and an enhanced optical-to-THz conversion efficiency. However, nanostructures also result in strong scattering of the optical EM wave, making the analytical expression used for carrier generation rate inaccurate. In addition, the geometrically complicated metallic nanostructures and the resulting spatially fast varying plasmon modes cannot be discretized and accounted for accurately using an orthogonal FDTD grid (unless an extremely fine grid, which would prohibitively increase the computational requirements, is used). 

Methods other than FDTD are also developed for simulating conventional PCDs~\cite{Neshat2010, Khabiri2012, Khiabani2013, Young2014}. In~\cite{Neshat2010, Khabiri2012}, the frequency-domain optical EM field is used to compute the (time-modulated) carrier generation and the transient carrier response is simulated using a coupled Poisson-DD solver. In~\cite{Khiabani2013}, an equivalent circuit model is used to predict the optical-to-THz conversion efficiency. In~\cite{Young2014}, the semiconductor device is treated as a circuit attached to a THz antenna, on which the THz EM wave radiation are accounted for using a discontinuous Galerkin time-domain (DGTD) scheme. Even though these approaches have been shown to be efficient for the simulation of conventional PCDs, they cannot be used for simulating plasmonic PCDs since the approximations involved in their formulation and implementation do not allow for accurate modeling of the optical EM wave/field interactions on the nanostructures. 

To simulate plasmonic PCDs, numerical frameworks that rely on the finite element method (FEM) have been developed~\cite{Burford2016, Bashirpour2017}. FEM, thanks to its capability to operate on unstructured meshes, can accurately account for the EM field interactions on the metallic nanostructures. These frameworks compute the space distribution of the carrier generation using the EM field distribution obtained by a frequency-domain FEM. The time dependency of the carrier generation is assumed to be a pulse with a closed-form analytical expression. Then the carrier generation (which is now assumed to be known in space and time) is fed into the time-domain simulation of carrier dynamics. Even though the EM field interactions on metallic nanostructures are accurately accounted for with this approach, the nonlinear two-way coupling between the EM interactions and carrier dynamics is not fully considered. Furthermore, the space-charge screening effect, i.e., the fact that separation of electrons and holes induces a polarization vector that cancels the electric field, which results from the time-dependent behavior of the carrier distributions, is not captured by this approach.

In this work, we propose a multiphysics framework that makes use of discontinuous Galerkin (DG) methods~\cite{Cockburn1999, Hesthaven2008, Shu2016High} to simulate plasmonic PCDs. This framework does not suffer from the shortcomings of the previously-developed methods as briefly described above. Following the two stages of the PCD operation, the semiconductor device is numerically modeled also in two stages: The stationary state is described by a coupled system of the (nonlinear) Poisson equation and the bipolar DD equations. This coupled system is solved iteratively using the Gummel method and the linearized set of equations at each iteration are discretized using (stationary) DG schemes~\cite{Cockburn1998, Castillo2000, Chen2020steadystate}. The transient stage is described by a nonlinearly coupled system of the time-dependent Maxwell and DD equations. This coupled system is discretized using DGTD schemes~\cite{Hesthaven2008, Hesthaven2002, Liu2007, Xu2010local}, where the nonlinear coupling between the two sets of equations is fully taken into account during the time integration. In this numerical framework, we prefer to use DG-based discretization because (i) it provides high flexibility in meshing (just like FEM), (ii) it allows for explicit time integration (just like FDTD), (iii) it permits easy implementation of high-order (spatial) basis functions and high-order time integration schemes, and (iv) non-conformal meshes, adaptive \textit{hp}-refinement, and local time stepping can be used to further reduce computational costs. These properties make DG methods very suitable for multiphysics and multiscale simulations~\cite{Hesthaven2008, Shu2016High, Jacobs2006, Harmon2016, Yan2016}. Explicit time integration, which is enabled by the use of DG-based discretization, equips the resulting multiphysics framework with several advantages:
(i) The nonlinear coupling between the Maxwell and DD equations can be easily and fully implemented. 
(ii) Each set of equations can use their own time-step size. Note that using a global time-step size would unnecessarily increase the computational requirements due to the large difference between the characteristic time scales of the Maxwell and DD equations. 
(iii) Since the explicit time integration does not require any matrix inversion, it can be parallelized easily and executed efficiently on distributed-memory computer clusters.

The rest of the paper is organized as follows. Section~\ref{model} describes the physical processes involved in the operation of the plasmonic PCDs and presents the systems of equations that are used to mathematically model them. Section~\ref{MuPDG} details the numerical schemes that are used for discretizing these systems and solving the resulting matrix equations. Some comments about the formulation and implementation are presented in Section~\ref{comments}. Section~\ref{results} demonstrates the accuracy, the efficiency, and the applicability of the proposed numerical framework via numerical examples. Finally, Section~\ref{conclusion} summarizes the work described in the paper and presents several future research directions.

\section{Formulation}
\label{sec:formulation}
\subsection{Physical and Mathematical Description}
\label{model}
The semiconductor most commonly used in PCD designs is the low-temperature grown gallium arsenide (LT-GaAs). LT-GaAs has an electron trapping time of less than $1~\mathrm{ps}$ and a direct optical band gap of 1.42eV. It can ``absorb'' optical EM energy and generate carriers with a lifetime shorter than $1~\mathrm{ps}$.
Two conductive electrodes are deposited on the semiconductor substrate. A bias voltage, which is applied to these electrodes, drives the carriers towards them, resulting in a THz current. Consequently, the electrodes act as a current source attached to a THz antenna. For a conventional PCD, the optical-to-THz conversion efficiency is limited by the amount of absorbed optical EM energy, which is usually small due to the high refractive index of the semiconductor layer.

Plasmonic PCDs utilize metallic nanostructures to increase the optical-to-THz conversion efficiency. These nanostructures are designed to introduce plasmon modes with resonances around the operating frequency of the optical EM wave excitation. Highly localized and strong EM fields associated with these modes significantly increase the carrier generation inside the semiconductor layer. In addition to introducing the plasmon modes, the nanostructures, which are either placed between the electrodes or etched on to them~\cite{Lepeshov2017review, Burford2017review, Kang2018review}, change the static electric field distribution. For the latter case, they act as a part of the electrode they are attached to and effectively reduce the distance between the carriers and the electrodes. Both of these effects significantly change the current~\cite{Moon2012,Yardimci2018review}. Therefore, the nanostructures are designed not only to generate strong localized EM fields, but also to maximize the overall current response. We also should mention here that the impedance mismatch between the semiconductor device (acting as the current source) and the THz antenna is not as critical as the current generation rate in terms of PCD’s performance~\cite{Lepeshov2017review, Yardimci2018review, Burford2016}.

Mathematical modeling of a PCD has to account for carrier dynamics, EM interactions, and their nonlinear coupling. While Maxwell equations are used to describe EM interactions, various models ranging from semi-classical to full-quantum approaches have been developed to describe carrier dynamics. For PCDs, the semi-classical DD model can be used accurately since the device size ($\sim 10 \mu$m) is far larger than the mean free path of the carriers ($\sim\! 0.01 \mu$m)~\cite{Chuang2012, Vasileska2010, Moreno20141}. In this work, we prefer to use the bipolar DD equations since they can account for the difference in the physical parameters associated with electrons and holes and allow for modeling the space-charge screening effect. Thus, carrier dynamics and EM interactions that happen on a PCD are mathematically described by the following coupled system~\cite{Vasileska2010, Chuang2012, Sha2012, In2014, Sha2012grating, Sha2014}
\begin{align}
\label{t_H0} & \mu ({\mathbf{r}}){\partial _t}{\mathbf{H}}({\mathbf{r}},t) =  - \nabla  \times {\mathbf{E}}({\mathbf{r}},t)\\
\label{t_E0} & \varepsilon ({\mathbf{r}}){\partial _t}{\mathbf{E}}({\mathbf{r}},t) = \nabla  \times {\mathbf{H}}({\mathbf{r}},t) - [{{\mathbf{J}}_e}({\mathbf{r}},t) + {{\mathbf{J}}_h}({\mathbf{r}},t)]\\
\label{t_N0} & q{\partial _t}{n_c}({\mathbf{r}},t) =  \pm \nabla  \cdot {{\mathbf{J}}_c}({\mathbf{r}},t) - q[R({n_e},{n_h}) - G({\mathbf{E}},{\mathbf{H}})]\\
\label{t_J0} & {{\mathbf{J}}_c}({\mathbf{r}},t) = q{\mu _c}({\mathbf{E}}){\mathbf{E}}({\mathbf{r}},t){n_c}({\mathbf{r}},t) \pm q{d_c}({\mathbf{E}})\nabla {n_c}({\mathbf{r}},t).
\end{align}
Here, ${\mathbf{r}}$ is the location vector, subscript $c\in \{ e,h\}$ represents the carrier type and hereinafter the upper and lower signs should be selected for electron ($c=e$) and hole ($c=h$), respectively, ${\mathbf{E}}({\mathbf{r}},t)$ and ${\mathbf{H}}({\mathbf{r}},t)$ are the electric and magnetic field intensities, $\varepsilon ({\mathbf{r}})$ and $\mu ({\mathbf{r}})$ are the dielectric permittivity and permeability, ${n_e}({\mathbf{r}},t)$ and ${n_h}({\mathbf{r}},t)$ are the electron and hole densities, ${{\mathbf{J}}_e}({\mathbf{r}},t)$ and ${{\mathbf{J}}_h}({\mathbf{r}},t)$ are the current densities due to electron and hole movement, $R({n_e},{n_h})$ and $G({\mathbf{E}},{\mathbf{H}})$ are the recombination and generation rates, ${\mu _e}({\mathbf{E}})$ and ${\mu _h}({\mathbf{E}})$ are the field-dependent electron and hole mobilities, ${d_e}({\mathbf{E}}) = {V_T}{\mu _e}({\mathbf{E}})$ and ${d_h}({\mathbf{E}}) = {V_T}{\mu _h}({\mathbf{E}})$ are the electron and hole diffusion coefficients, respectively. {\color{black} Here, $V_T=k_B T/q$ is the thermal voltage, $k_B$ is the Boltzmann constant, and $T$ is the absolute temperature.}

In this work, we use the trap assisted recombination described by the Shockley-Read-Hall model~\cite{Chuang2012, Vasileska2010, Chen2020steadystate}
{\color{black}
\begin{equation}
\label{recomb}
{R}({n_e},{n_h}) =\frac{ {n_e}({\mathbf{r}}){n_h}({\mathbf{r}}) - {n_i^2} } { {\tau _e}[{n_{h1}} + {n_h}({\mathbf{r}})] + {\tau _h}[{n_{e1}} + {n_e}({\mathbf{r}})] }
\end{equation}
where ${n_i}$ is the intrinsic carrier concentration, ${\tau _e}$ and ${\tau _h}$ are the carrier lifetimes, and ${n_{e1}}$ and ${n_{h1}}$ are model parameters related to the trap energy level. The generation rate $G({\mathbf{E}},{\mathbf{H}})$ is a function of the optical EM fields and is introduced below in~\eqref{generation}. The parallel-field dependent mobility model~\cite{Moreno20141} is used to account for the field-dependency and velocity-saturation for a more accurate prediction of the generated THz current (mobility affects the carrier drift velocity, which in return has a big impact on the current).}

The bias voltage is applied to the electrodes continuously throughout the operation of the device. Before the optical EM wave excitation is turned on, the device is under a non-equilibrium stationary state~\cite{Chen2020steadystate}. When the EM wave impinges on the device, carriers are generated. However, the density of these carriers is several orders of magnitude smaller than the doping concentration and the stationary state of the device is assumed to be only weakly perturbed by the optical EM wave excitation. Therefore, field intensities and carrier and current densities are separated into stationary and transient components as 
${\mathbf{E}}({\mathbf{r}},t) = {{\mathbf{E}}^s}({\mathbf{r}}) + {{\mathbf{E}}^t}({\mathbf{r}},t)$,
${\mathbf{H}}({\mathbf{r}},t) = {{\mathbf{H}}^s}({\mathbf{r}}) + {{\mathbf{H}}^t}({\mathbf{r}},t)$,
${n_c}({\mathbf{r}},t) = n_c^s({\mathbf{r}}) + n_c^t({\mathbf{r}},t)$, ${{\mathbf{J}}_c}({\mathbf{r}},t) = {\mathbf{J}}_c^s({\mathbf{r}}) + {\mathbf{J}}_c^t({\mathbf{r}},t)$, respectively. Here, the superscript ``\textit{s}'' and ``\textit{t}'' stands for stationary and transient components, respectively.
Similarly, the recombination rate in \eqref{t_N0} is decomposed into stationary and transient components as $R({n_e},{n_h}) = {R^s}(n_e^s,n_h^s) + {R^t}(n_e^t,n_h^t)$~\cite{Moreno2014}. Since the static electric field component is at least two orders of magnitude stronger than its transient counterpart, $\left| {{\mathbf{E}}^{s}}(\mathbf{r}) \right|\gg \left| {{\mathbf{E}}^{t}}(\mathbf{r},t) \right|$, it is assumed that mobility is only a function of ${{\mathbf{E}}^{s}}(\mathbf{r})$, i.e., ${{\mu }_{c}}(\mathbf{E})\approx {{\mu }_{c}}({{\mathbf{E}}^{s}})$, $c\in \{e,h\}$.

The generation rate in~\eqref{t_N0} is defined as~\cite{Chuang2012}
\begin{equation}
\label{generation}
G({{\mathbf{E}}^t},{{\mathbf{H}}^t}) = \frac{\eta \alpha \lambda}{hc} \left| {{{\mathbf{E}}^t}({\mathbf{r}},t) \times {{\mathbf{H}}^t}({\mathbf{r}},t)} \right|
\end{equation}
where $\lambda$ is the wavelength of the optical EM wave excitation, $\alpha $ is the photon absorption coefficient, $\eta$ is the quantum yield, $h$ is the Planck constant, and $c$ is the light speed. We should emphasize here that $G({{\mathbf{E}}^t},{{\mathbf{H}}^t})$ depends only on the transient EM fields ${{\mathbf{E}}^t({\mathbf{r}},t)}$ and ${{\mathbf{H}}^t({\mathbf{r}},t)}$ that are generated on the PCD due to the optical EM wave excitation and include the effect of all transient EM interactions.


Following the discussion above, the coupled of system of time-dependent Maxwell and DD equations in~\eqref{t_H0}-\eqref{t_J0} is effectively decomposed into two coupled sets of equations. The stationary components of the variables (the ones with superscript ``\textit{s}'') satisfy a coupled system of Poisson and stationary DD equations. This is the first equation system and it reads
\begin{align}
\label{s_PS} 
& \nabla  \cdot [\varepsilon ({\mathbf{r}}) {{\mathbf{E}}^s}({\mathbf{r}}) ] = q[C({\mathbf{r}}) + n_h^s({\mathbf{r}}) - n_e^s({\mathbf{r}})] \\
\label{s_C}  
& \nabla  \cdot {\mathbf{J}}_c^s({\mathbf{r}}) =  \pm q{R^s}(n_e^s,n_h^s) \\
\label{s_J}  
& {\mathbf{J}}_c^s({\mathbf{r}}) = q{\mu _c}({{\mathbf{E}}^s}){{\mathbf{E}}^s}({\mathbf{r}})n_c^s({\mathbf{r}}) \pm q{d_c}({{\mathbf{E}}^s})\nabla n_c^s({\mathbf{r}})
\end{align}
{\color{black} where $C({\mathbf{r}})$ is the doping concentration.}
Eliminating the stationary components in~\eqref{t_H0}-\eqref{t_J0} yields a reduced coupled system of Maxwell and DD equations. This is the second equation system and it reads
\begin{align}
\label{t_H} & \mu(\mathbf{r}) {\partial _t}{{\mathbf{H}}^t}({\mathbf{r}},t) = -\nabla \times {{\mathbf{E}}^t} ({\mathbf{r}},t) \\
\label{t_E} & \varepsilon(\mathbf{r}) {\partial _t}{{\mathbf{E}}^t}({\mathbf{r}},t) \!=\!  \nabla \times {{\mathbf{H}}^t}({\mathbf{r}},t) \!-\! [{\mathbf{J}}_e^t({\mathbf{r}},t) + {\mathbf{J}}_h^t({\mathbf{r}},t)]\\
\label{t_N} & q{\partial _t}n_c^t({\mathbf{r}},t) \!=\! \pm \nabla \cdot {\mathbf{J}}_c^t({\mathbf{r}},t) \!-\! q[{R^t}(n_e^t,n_h^t) \!-\! G({{\mathbf{E}}^t},{{\mathbf{H}}^t})]\\
\label{t_J} & {\mathbf{J}}_c^t({\mathbf{r}},t) = q{\mu _c}({\mathbf{E}}^s)([{{\mathbf{E}}^s}({\mathbf{r}}) + {{\mathbf{E}}^t}({\mathbf{r}},t)]n_c^t({\mathbf{r}},t) + {{\mathbf{E}}^t}({\mathbf{r}},t) n_c^s({\mathbf{r}},t)) \pm q{d_c}({\mathbf{E}}^s)\nabla n_c^t({\mathbf{r}},t).
\end{align}
Equations \eqref{t_H}-\eqref{t_J} are a strongly-coupled nonlinear system. In \eqref{t_N}, $G({{\mathbf{E}}^t},{{\mathbf{H}}^t})$ leads to an exponential increase in carrier densities. In \eqref{t_J}, carriers are driven by both ${{\mathbf{E}}^t}({\mathbf{r}},t)$ and ${{\mathbf{E}}^s}({\mathbf{r}})$. The carrier motion produces the free current densities ${\mathbf{J}}_e^t({\mathbf{r}},t)$ and ${\mathbf{J}}_h^t({\mathbf{r}},t)$, which contributes to the EM wave/field interactions described by~\eqref{t_H}-\eqref{t_E}. In return, these EM wave/field interactions, change $G({{\mathbf{E}}^t},{{\mathbf{H}}^t})$. We should also note here that ${R^t}(n_e^t,n_h^t)$ is also a nonlinear function, which becomes significant when carrier densities are high, balancing the carrier generation.

\subsection{The Multiphysics DG Solver}
\label{MuPDG}
A complete simulator for numerical characterization of PCDs consists of a stationary scheme that solves  the Poisson-DD system \eqref{s_PS}-\eqref{s_J} and a time-domain scheme that solves the Maxwell-DD system \eqref{t_H}-\eqref{t_J}. The stationary solutions are used as inputs to the time-domain scheme. The stationary scheme uses the Gummel iteration method to take the nonlinearity into account. The linearized set of equations that needs to be solved at every iteration of this method is discretized using a stationary DG scheme. For details, we refer the reader to~\cite{Chen2020steadystate}. In the rest of this paper, we focus on the time-domain scheme proposed to solve the Maxwell-DD system \eqref{t_H}-\eqref{t_J}.

\subsubsection{Nonlinear Coupling}
\label{coupling}
To integrate the time-dependent Maxwell-DD system in time, we use an explicit scheme. To account for the difference in the characteristic time scales of the EM interactions and the carrier dynamics in an efficient manner, the Maxwell and DD equations are updated using two different time-step sizes. More specifically, since the carrier response is much slower than the time variation of the EM waves/fields, the DD equations are updated using a larger time-step size. The nonlinear coupling between the two equation sets is accounted for by alternately feeding the updated solutions into each other.

\begin{figure}[!t]
	\centerline{\includegraphics[width=0.6\columnwidth]
		{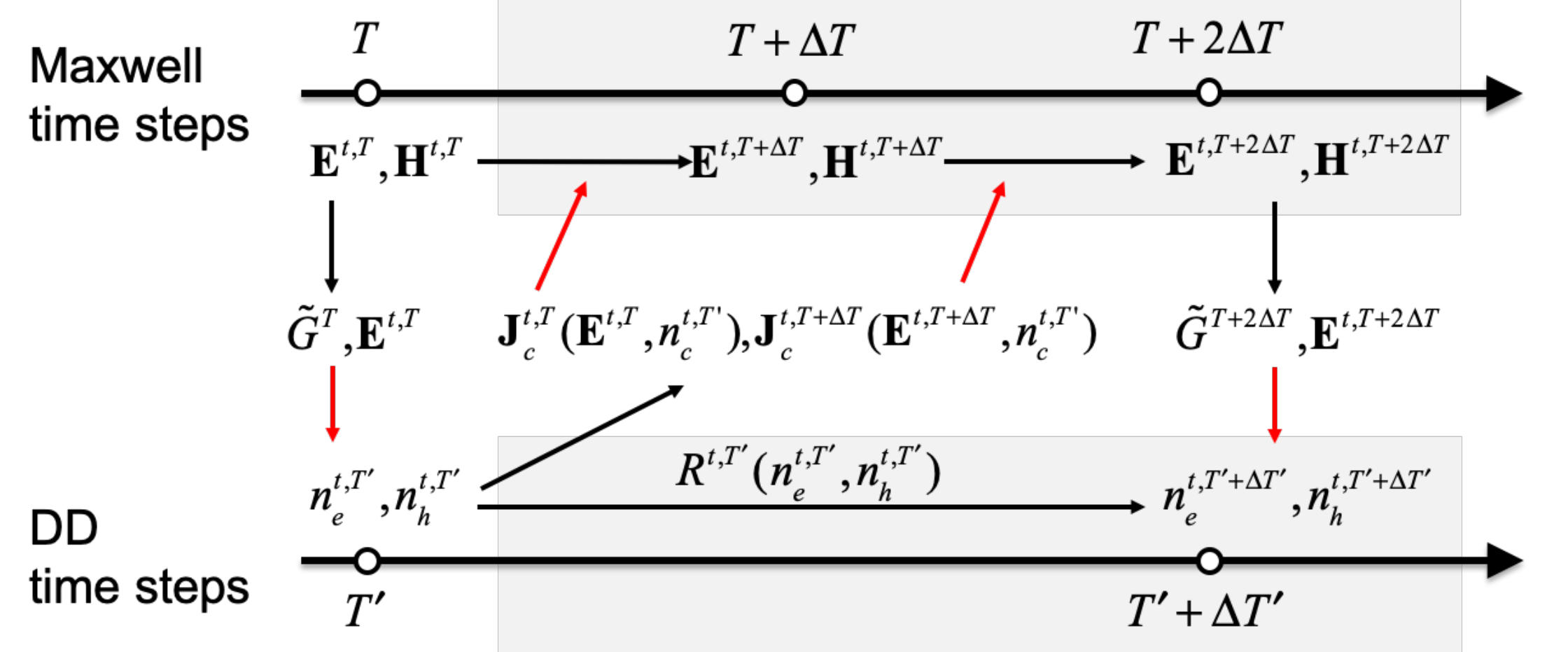}}
	\caption{Time integration of the Maxwell-DD system.}
	\label{TimeIntegration}
\end{figure}

This time marching scheme is shown in Fig.~\ref{TimeIntegration}, where the time-step size for the DD equations is assumed to be twice the step size for the Maxwell equations for illustration purposes. Let us suppose that two equation systems are updated separately using two different time integration schemes (to be discussed in the next section) and denote the time steps and step sizes of the Maxwell and DD systems with \{$T$, $\Delta T$\} and \{$T'$, $\Delta T'$\}, respectively. We should note here that, in what follows, subscripts $T$ and $T'$ mean that the variables they are attached to are computed/sampled/updated at $t=T$ and $t=T'$, respectively. 

Let us assume that the time steps of the two systems are synchronized at time $T=T'$. {\color{black} First, EM fields $\{\mathbf{E}^{t,T-\Delta T},\mathbf{H}^{t,T-\Delta T}\}$ and $\{\mathbf{E}^{t,T}, \mathbf{H}^{t,T}\}$ are used in~\eqref{generation} to compute the generation rates $G^{T-\Delta T}$ and ${G^T}$, respectively. The generation rate is averaged over times $T-\Delta T$ and $T$, i.e., $\tilde G^T=(G^{T-\Delta T}+G^T)/2$. Carrier densities $n_c^{t,T'-\Delta T'}$ are used in~\eqref{recomb} to compute the recombination rate $R^{t,T'-\Delta T'}$. Then, $\tilde G^T$ and $R^{t,T'-\Delta T}$ are used together in~\eqref{t_N} to update $n_c^{t,T'}$ (from step $T'-\Delta T'$ to $T'$).} Then, $n_c^{t,T'}$ are used in~\eqref{t_J} to compute  the current densities $\mathbf{J}_c^{t,T}$ and $\mathbf{J}_c^{t,T+\Delta T}$. $\mathbf{J}_c^{t,T}$ and $\mathbf{J}_c^{t,T+\Delta T}$ are used in~\eqref{t_E} to update the Maxwell equations at steps $T$ and $T+\Delta T$ to produce $\{\mathbf{E}^{t,T+\Delta T },\mathbf{H}^{t,T+\Delta T}\}$ and $\{\mathbf{E}^{t,T+2\Delta T}, \mathbf{H}^{t, T+2\Delta T}\}$, respectively. The time steps of the two systems match again at time $T+2\Delta T = T'+\Delta T'$, {\color{black} $\tilde G^{T+2\Delta T}=(G^{T+\Delta T}+G^{T+2\Delta T})/2$ and $R^{t,T'}$ are computed using $\{\mathbf{E}^{t,T+\Delta T },\mathbf{H}^{t,T+\Delta T}\}$ and $\{\mathbf{E}^{t,T+2\Delta T}, \mathbf{H}^{t, T+2\Delta T}\}$ and $n_c^{t,T'}$, respectively,} and the process described here is repeated. 

\subsubsection{Discretization}
\label{discretization}
The time-dependent DD and Maxwell equations [\eqref{t_N}-\eqref{t_J} and~\eqref{t_H}-\eqref{t_E}, respectively] are discretized using DG schemes in time-domain. We start the description of discretization with the DD equations~\eqref{t_N}-\eqref{t_J}. Since the electron and hole DD equations only differ by the sign in front of the drift term, we only discuss the electron DD equation. Also we note that we drop the subscript ``\textit{e}'' (meaning electron) and the superscript ``\textit{t}'' (meaning transient) from the variables to simplify the notation. Since $\mathbf{E}^s(\mathbf{r})$ is provided as an input to the time-domain simulation, variables that are functions of it are assumed to change with $\mathbf{r}$ only. Furthermore, as the three drift terms in~\eqref{t_J} are treated in the same way, for brevity, 
we combine them under one term that is denoted by ${\mathbf{v}}({\mathbf{r}},t)n({\mathbf{r}},t)$. Under those notation simplifications, the electron DD equations~\eqref{t_N}-\eqref{t_J} are expressed as the following initial-boundary value problem (IBVP)
\begin{align}
\label{IVP_DD0}
& {\partial _t}n({\mathbf{r}},t) = \nabla \cdot [d({\mathbf{r}}){\mathbf{q}}({\mathbf{r}},t)] + \nabla \cdot [{\mathbf{v}}({\mathbf{r}},t)n({\mathbf{r}},t)] - R({\mathbf{r}},t),{\mathbf{r}} \in \Omega \\
\label{IVP_DD1}
& {\mathbf{q}}({\mathbf{r}},t) = \nabla n({\mathbf{r}},t), {\mathbf{r}} \in \Omega\\
\label{IVP_DD2}
& n({\mathbf{r}},t) = {f_D} ({\mathbf{r}}), {\mathbf{r}} \in \partial {\Omega _D}\\
\label{IVP_DD3}
& \mathbf{\hat{n}}({\mathbf{r}}) \cdot [d({\mathbf{r}}) {\mathbf{q}} ({\mathbf{r}},t) + {\mathbf{v}}({\mathbf{r}},t)n({\mathbf{r}},t)] = {f_R}({\mathbf{r}}), {\mathbf{r}} \in \partial {\Omega _R}.
\end{align}
Here, ${\mathbf{q}}({\mathbf{r}},t)$ is an auxiliary variable introduced to reduce the order of the spatial derivative in the diffusion term, $R({\mathbf{r}},t) \equiv {R^t}(n_e^t,n_h^t) - G({{\mathbf{E}}^t},{{\mathbf{H}}^t})$, $\partial {\Omega _D}$ and $\partial {\Omega _R}$ represent the surfaces where Dirichlet and Robin boundary conditions are enforced and ${f_D}$ and ${f_R}$ are the coefficients associated with these boundary conditions, respectively, and ${\mathbf{\hat n}}({\mathbf{r}})$ denotes the outward normal vector $\partial {\Omega _R}$. For the problems considered in this work, $\partial {\Omega _D}$ represents the electrode/semiconductor interfaces and ${f_D}=0$; and $\partial {\Omega _R}$ represents the semiconductor/insulator interfaces and ${f_R}=0$ indicating no carrier spills out those interfaces~\cite{Schroeder1994}.

To facilitate the numerical solution of the IBVP described by \eqref{IVP_DD0}-\eqref{IVP_DD3}, $\Omega$ is discretized into $K$ non-overlapping tetrahedrons. The volumetric support of each of these discretization elements is represented by ${\Omega _k}$, $k = 1, \ldots ,K$. Let $\partial {\Omega _k}$ and ${\mathbf{\hat n}}({\mathbf{r}})$ denote the surface of ${\Omega _k}$ and the outward unit vector normal to $\partial {\Omega _k}$, respectively. Hereinafter, the nodal expansion~\cite{Hesthaven2008} is used. Testing \eqref{IVP_DD0} and \eqref{IVP_DD1} with Lagrange polynomials, ${\ell _i}({\mathbf{r}})$, $i = 1, \ldots ,{N_p}$~\cite{Hesthaven2008}, on element $k$ and applying the divergence theorem yield the weak form
\begin{align}
\nonumber & \int_{{\Omega _k}} {\partial _t} {n_k({\mathbf{r}},t){\ell _i}({\mathbf{r}})dV} = - \int_{{\Omega _k}} {R({\mathbf{r}},t){\ell _i}({\mathbf{r}})dV} 
-\! \int_{{\Omega _k}}\!\!\!\! {d({\mathbf{r}}){\mathbf{q}}_k({\mathbf{r}},t) \! \cdot \! \nabla {\ell _i}({\mathbf{r}})dV}  
\! + \! \oint_{\partial {\Omega _k}}\!\!\!\! {{\mathbf{\hat n}}({\mathbf{r}}) \! \cdot \! {{  (  d{\mathbf{q}}  )  }^*}{\ell _i}({\mathbf{r}})dS} \\
\label{weakN} & \quad \quad \quad \quad \quad \quad \quad \quad \quad \quad -\! \int_{{\Omega _k}}\!\!\!\! {{\mathbf{v}}({\mathbf{r}},t){n_k}({\mathbf{r}},t) \! \cdot \! \nabla {\ell _i}({\mathbf{r}})dV} 
\! + \! \oint_{\partial {\Omega _k}}\!\!\!\! {{\mathbf{\hat n}}({\mathbf{r}}) \! \cdot \! {{( {\mathbf{v}} n )}^*}{\ell _i}({\mathbf{r}})dS} \\
\label{weakQ} & \int_{{\Omega _k}} {{q_k^\nu} } { ({\mathbf{r}},t) {\ell _i}({\mathbf{r}})dV} = - \int_{{\Omega _k}} {{n_k}({\mathbf{r}},t)  {\partial _{\nu} }  {\ell _i}({\mathbf{r}})dV} 
+ \! \oint_{\partial {\Omega _k}} {{{\hat n}_\nu }({\mathbf{r}})n^*{\ell _i}({\mathbf{r}})dS}.
\end{align}
Here, ${N_p} = (p + 1)(p + 2)(p + 3)/6$ is the number of interpolation nodes, $p$ is the order of the Lagrange polynomials and $\nu \in \{x,y,z\}$ is used for identifying the components of the vectors in the Cartesian coordinate system. We note here ${n_k}({\mathbf{r}},t)$ and ${\mathbf{q}_k}({\mathbf{r}},t)$ denote the local solutions on element $k$ and the global solutions on $\Omega$  are the direct sum of the local solutions.

In~\eqref{weakN} and~\eqref{weakQ}, ${n^*}$, ${(d{\mathbf{q}})^*}$, and ${({\mathbf{v}}n)^*}$ are numerical fluxes ``connecting'' element $k$ to its neighboring elements. The variables involved in the definition of the numerical fluxes reside on the interfaces between elements and their explicit dependencies on $k$,  ${\mathbf{r}}$, and $t$ are dropped to simplify the notation. For the diffusion term, the local DG (LDG) alternate flux~\cite{Cockburn1998} is used for ${n^*}$ and ${(d{\mathbf{q}})^*}$. For the drift term, the local Lax-Friedrichs flux~\cite{Hesthaven2008} is used for ${({\mathbf{v}}n)^*}$. On the surfaces of the element,
they are defined as
\begin{align*}
& {n^*} = \left\{ n \right\} + 0.5\boldsymbol{\hat \beta} \cdot {\mathbf{\hat n}} \left[\kern-0.15em\left[ n 
\right]\kern-0.15em\right] \\
& {\left( {d{\mathbf{q}}} \right)^*} = \left\{ {d{\mathbf{q}} } \right\} - 0.5\boldsymbol{\hat \beta} ({\mathbf{\hat n}} \cdot \left[\kern-0.15em\left[ {d{\mathbf{q}}} \right]\kern-0.15em\right]) \\
& {\left( {{\mathbf{v}}n} \right)^*} = \left\{ {{\mathbf{v}}n} \right\} + \alpha {\mathbf{\hat n}}\left[\kern-0.15em\left[ n 
\right]\kern-0.15em\right].
\end{align*}
Here, ``average'' and ``jump'' operators are respectively defined as $\left\{\odot \right\} = 0.5({\odot ^ - } + { \odot ^ + })$ and $\left[\kern-0.15em\left[ \odot 
\right]\kern-0.15em\right] = { \odot ^ - } - { \odot ^ + }$, where $ \odot $ is a scalar or a vector variable. Superscripts ``-'' and ``+'' refer to variables defined in element $k$ and in its neighboring element, respectively. The vector $\boldsymbol{\hat \beta}$ determines the upwind direction of $n$ and $(d\mathbf{q})$. In LDG, it is essential to choose opposite directions for $n$ and $(d\mathbf{q})$ (note the sign difference in the definitions above), while the precise direction of each variable is not important~\cite{Cockburn1998, Hesthaven2008, Shu2016}. In this work, we choose $\boldsymbol{\hat \beta}={\mathbf{\hat n}}$ on each element surface. The local Lax-Friedrichs, with $\alpha  = \max (|{\mathbf{\hat n}} \cdot {{\mathbf{v}}^ - }|,|{\mathbf{\hat n}} \cdot {{\mathbf{v}}^ + }|)/2$~\cite{Hesthaven2008}, mimics the path of the information propagation. On $\partial {\Omega _D}$, ${n^*} = {f_D}$, ${(d{\mathbf{q}})^*} = {(d{\mathbf{q}})^ - }$ and ${({\mathbf{v}}n)^*} = {\mathbf{v}^-}{f_D}$, and on $\partial {\Omega _R}$, ${n^*} = {n^ - }$ and ${(d{\mathbf{q}})^*} + {({\mathbf{v}}n)^*} = {f_R}$. We note that ${(d{\mathbf{q}})^*}$ and ${({\mathbf{v}}n)^*}$ are not assigned independently on $\partial {\Omega _R}$.

${n_k}({\mathbf{r}},t)$ and ${q_k^{\nu}}({\mathbf{r}},t)$ are expanded using Lagrange polynomials ${\ell _i}({\mathbf{r}})$~\cite{Hesthaven2008}
\begin{align}
\label{expN}
& {n_k}({\mathbf{r}},t) \simeq \sum\limits_{i = 1}^{{N_p}} {n_k({{\mathbf{r}}_i},t){\ell _i}({\mathbf{r}})}  = \sum\limits_{i = 1}^{{N_p}} {n_{k,i}(t){\ell _i}({\mathbf{r}})}\\
\label{expQ}
& q_k^\nu ({\mathbf{r}},t) \simeq \mathop \sum \limits_{i = 1}^{{N_p}} {q_k^\nu }({{\mathbf{r}}_i},t){\ell _i}({\mathbf{r}}) = \mathop \sum \limits_{i = 1}^{{N_p}} q_{k,i}^{\nu}(t){\ell _i}({\mathbf{r}})
\end{align}
where $\mathbf{r}_i$, $i=1,\ldots,N_p$, denote the locations of the interpolation nodes, $n_{k,i}(t)$ and $q_{k,i}^{\nu}(t)$, $\nu  \in \{ x,y,z\} $, $k = 1, \ldots ,K$, are the unknown coefficients to be solved for. Substituting \eqref{expN} and \eqref{expQ} into \eqref{weakN} and \eqref{weakQ}, we obtain the semi-discretized form
\begin{align}
\label{semiDD0} &
{{\bar M}_{k}}{\partial _t}{{\bar N}_k(t)} = {{\bar C}_{k}}{{\bar N}_k(t)} + {{\bar C}_{kk'}}{{\bar N}_{k'}(t)} + { {\bar D}_k {\bar d}_k }{{\bar Q}_k(t)} + { {\bar D}_{kk'} {\bar d}_{k'} }{{\bar Q}_{k'}(t)} - \bar B_k^n(t) \\
\label{semiDD1} &
\bar M_{k}^{\mathbf{q}}{{\bar Q}_k(t)} = {{\bar G}_{k}}{{\bar N}_k}(t) + {{\bar G}_{kk'}}{{\bar N}_{k'}}(t) + \bar B_k^{\mathbf{q}}
\end{align}
where the unknown vectors are defined as ${\bar N_k(t)} = [n_{k,1}(t),...,n_{k,N_p}(t)]^T$ and ${\bar Q}_k(t) = {[\bar Q_k^x(t),\bar Q_k^y(t), \bar Q_k^z(t)]^T}$, with $\bar Q_k^\nu(t)  = [q_{k,1}^{\nu}(t),...,q_{k,{N_p}}^{\nu}(t)]$, $\nu  \in \{ x,y,z\} $.

In~\eqref{semiDD0}-\eqref{semiDD1}, ${\bar M}_k$ and ${{\bar M}_k^{\mathbf{q}}}$ are mass matrices. $\bar M_{k}^{\mathbf{q}}$ is a $3 \times 3$ block diagonal matrix with blocks ${\bar M_{k}}$. ${\bar M_{k}}$ is defined as
\begin{equation*}
\bar M_{k}(i,j) = \int_{{\Omega _k}} {{\ell _i}({\mathbf{r}}){\ell _j}({\mathbf{r}})} dV, \quad i,j = 1, \ldots, {N_p}.
\end{equation*}
$\bar{d}_k$ is a diagonal matrix with entries ${d_1},...,{d_K}$, where ${d_k} = (d_k^x,d_k^y,d_k^z)$ and $d_k^\nu (i) = {d_k}({{\mathbf{r}}_i})$, $i = 1,...,{N_p}$, $\nu \in \{x,y,z\}$. We note that $d({\mathbf{r}})$ is assumed isotropic and constant in each element.
Matrices ${\bar G}_k$ and ${\bar G}_{kk'}$, and ${\bar D}_k$ and ${\bar D}_{kk'}$ correspond to the gradient and the divergence operators, respectively.
For the LDG flux, ${\bar D}_k = -{{\bar G}_k^T}$ and ${\bar D}_{kk'} = -{{\bar G}_{kk'}^T}$. ${\bar G}_k$ is a $3{N_p} \times {N_p}$ matrix and it has contributions from the volume and surface integral terms on the right hand side of~\eqref{weakQ}: ${\bar G}_k=G_{k}^{\mathrm{V}} + G_{k}^{\mathrm{S}}$, $\bar G_{k}^{\mathrm{V}} = {\left[ {\bar S_k^x\;\bar S_k^y\;\bar S_k^z} \right]^T}$, $\bar G_{k}^{\mathrm{S}} = {\left[ {\bar L_k^x\;\bar L_k^y\;\bar L_k^z} \right]^T}$, where
\begin{equation*}
\bar S_k^\nu (i,j) =  - \int_{{\Omega _k}} { \partial_\nu {{\ell _i}({\mathbf{r}})}  {\ell _j}({\mathbf{r}})} dV
\end{equation*}
and
\begin{equation*}
\bar L_k^\nu (i,j) = \frac{{1 + sign(\boldsymbol{\hat \beta} \cdot {\mathbf{\hat n}})}}{2}{\theta _k}(j)\oint_{\partial {\Omega _{kk'}}} \!\!\!\! {{{\hat n}_\nu }({\mathbf{r}}){\ell _i}({\mathbf{r}}){\ell _j}({\mathbf{r}})dS}.
\end{equation*}
Here, $\partial {\Omega _{kk'}}$ denotes the interface connecting element $k$ and $k'$ and ${\theta _k}(j)$ selects the interpolation nodes on the interface
\begin{equation*}
{\theta _k}(j) = \left\{ \begin{gathered}
1,\quad {{\mathbf{r}}_j} \in {\Omega _k},{{\mathbf{r}}_j} \in \partial {\Omega _{kk'}} \hfill \\
0,\quad \mathrm{otherwise} \hfill \\ 
\end{gathered}  \right. .
\end{equation*}
Matrix ${\bar G_{kk'}}$ corresponds to the surface integral term in \eqref{weakQ}, which involves the unknowns from neighboring elements of element $k$: $\bar G_{kk'} = {\left[ {\bar L_{k'}^x\;\bar L_{k'}^y\;\bar L_{k'}^z} \right]^T}$, where
\begin{equation*}
\bar L_{k'}^\nu (i,j) = \frac{{1 - sign(\boldsymbol{\hat \beta} \cdot {\mathbf{\hat n}})}}{2}{\theta _{k'}}(j)\oint_{\partial {\Omega _{kk'}}} \!\!\!\! {{{\hat n}_\nu }({\mathbf{r}}){\ell _i}({\mathbf{r}}){\ell _j}({\mathbf{r}})dS}.
\end{equation*}
Similarly, matrix ${\bar C}_k$ has contributions from the fourth term (the volume integral) and the fifth term (the surface integral) on the right hand side of~\eqref{weakN}: ${\bar C}_k = \bar C_{k}^{\mathrm{V}} + \bar C_{k}^{\mathrm{S}}$, where
\begin{align*}
\bar C_{k}^{\mathrm{V}}(i,j) = - \sum\nolimits_\nu \int_{{\Omega _k}} {{v_\nu }({\mathbf{r}},t)  {\partial_\nu {\ell _i}({\mathbf{r}})} {\ell _j}({\mathbf{r}})} dV
\end{align*}
and
\begin{align*}
\nonumber \bar C_{k}^{\mathrm{S}} & (i,j) = {\theta _k}(j)
\oint_{\partial {\Omega _{kk'}}} {[ \frac{1}{2}\sum\nolimits_\nu {{{\hat n}_\nu }({\mathbf{r}}){v_\nu }({\mathbf{r}},t) + \alpha ({\mathbf{r}},t)} ] {\ell _i}({\mathbf{r}}){\ell _j}({\mathbf{r}})dS} .
\end{align*}
Matrix ${\bar C_{kk'}}$ corresponds to the last surface integral term in \eqref{weakN}, which involves the unknowns from neighboring elements of the element $k$:
\begin{align*}
\nonumber \bar C_{kk'} & (i,j) = {\theta _{k'}}(j)
\oint_{\partial {\Omega _{kk'}}} {[ \frac{1}{2}\sum\nolimits_\nu  {{{\hat n}_\nu }({\mathbf{r}}){v_\nu }({\mathbf{r}},t) - \alpha ({\mathbf{r}},t)} ]{\ell _i}({\mathbf{r}}){\ell _j}({\mathbf{r}})dS} .
\end{align*}
Matrices $\bar B_k^{n}$ and $\bar B_k^{\mathbf{q}}$ are contributed from the force term and boundary conditions (where element $k'$ does not exist) and are expressed as
\begin{align*}
\nonumber & \bar B_k^n(i) = \int_{{\Omega _k}} {R({\mathbf{r}},t){\ell _i}({\mathbf{r}})dV} + \oint_{\partial {\Omega _k} \cap \partial {\Omega _R}} {{f_R}({\mathbf{r}}){\ell _i}({\mathbf{r}})dS} 
 + \oint_{\partial {\Omega _k} \cap \partial {\Omega _D}} {{\mathbf{\hat n}}({\mathbf{r}}) \cdot {\mathbf{v}}({\mathbf{r}},t){f_D}({\mathbf{r}}){\ell _i}({\mathbf{r}})dS}\\
& \bar B_k^{{\mathbf{q}},\nu }(i) = \oint_{\partial {\Omega _k} \cap \partial {\Omega _D}} {{{\hat n}_\nu }({\mathbf{r}}){f_D}({\mathbf{r}}){\ell _i}({\mathbf{r}})dS} .
\end{align*}

To integrate~\eqref{semiDD0}-\eqref{semiDD1} in time, an explicit third-order total-variation-diminishing (TVD) Runge-Kutta method~\cite{Shu1988} is used. The high-order accuracy of this scheme matches that of the spatial discretization. With initial value ${n}({\mathbf{r}},t=0)=0$, time samples of the unknown vector ${\bar N_k(t)}$ are obtained step by step in time.

The time-dependent Maxwell equations \eqref{t_H} and \eqref{t_E} are discretized using the nodal DG method~\cite{Hesthaven2008,Sirenko2012}. First the simulation domain is divided into $K$ non-overlapping tetrahedrons with volumetric support $\Omega_k$, $k=1,...,K$. Just like before, $\partial {\Omega _k}$ and ${\mathbf{\hat n}}({\mathbf{r}})$ denote the surface of ${\Omega _k}$ and the outward unit vector normal to $\partial {\Omega _k}$, respectively. Testing \eqref{t_H} and \eqref{t_E} with the Lagrange polynomials ${\ell _i}({\mathbf{r}})$, $i = 1, \ldots ,{N_p}$, on element $k$ and applying the divergence theorem twice yield the strong form~\cite{Hesthaven2008}
\begin{align}
\label{strongE} & \int_{{\Omega _k}} \!\! \varepsilon ({\mathbf{r}}){\partial _t} {\mathbf{E}_k}({\mathbf{r}},t) {\ell _i}({\mathbf{r}})dV = \int_{{\Omega _k}} \!\! \nabla  \times {\mathbf{H}_k}({\mathbf{r}},t){\ell _i}({\mathbf{r}})dV 
 \! - \!\! \oint_{\partial {\Omega _k}} \!\!\!\! {\mathbf{\hat n}}({\mathbf{r}}) \! \times \! [{\mathbf{H}_k}({\mathbf{r}},t) \! - \! {\mathbf{H}}^*]{\ell _i}({\mathbf{r}})dS \! - \!\! \int_{{\Omega _k}} \!\!\!  {{\mathbf{J}}_k}({\mathbf{r}},t){\ell _i}({\mathbf{r}})dV \\
\label{strongH} & \int_{{\Omega _k}} \!\! \mu ({\mathbf{r}}){\partial _t} {\mathbf{H}_k}({\mathbf{r}},t){\ell _i}({\mathbf{r}})dV =  - \int_{{\Omega _k}} \!\! \nabla  \times {\mathbf{E}_k}({\mathbf{r}},t){\ell _i}({\mathbf{r}})dV 
 \! + \!\! \oint_{\partial {\Omega _k}} \!\!\!\! {\mathbf{\hat n}}({\mathbf{r}}) \! \times \! [{\mathbf{E}_k}({\mathbf{r}},t) - {\mathbf{E}}^*]{\ell _i}({\mathbf{r}})dS.
\end{align}
Here, $\mathbf{E}_k({\mathbf{r}},t)$ and $\mathbf{H}_k({\mathbf{r}},t)$ are the local solutions on element $k$, and ${\mathbf{E}}^*$ and ${\mathbf{H}}^*$ are the upwind numerical fluxes ``connecting'' element $k$ to its neighboring elements. They are expressed as~\cite{Hesthaven2008}
\begin{align*}
&  {{\mathbf{E}}^*} = (2\left\{ {Y{\mathbf{E}}} \right\} - {\mathbf{\hat n}} \times \left[\kern-0.15em\left[ {\mathbf{H}} 
 \right]\kern-0.15em\right])/(2\left\{ Y \right\})  \\
&  {{\mathbf{H}}^*} = (2\left\{ {Z{\mathbf{H}}} \right\} + {\mathbf{\hat n}} \times \left[\kern-0.15em\left[ {\mathbf{E}} 
 \right]\kern-0.15em\right])/(2\left\{ Z \right\})
\end{align*}
where the average and the jump operators are the same as those defined before, $Z$ and $Y$ are the wave impedance and wave admittance, respectively. Same as before, the variables are defined on element surfaces and their explicit dependencies on $k$, ${\mathbf{r}}$, and $t$ are dropped to simplify the notation. On boundaries where there are no neighboring elements, the numerical fluxes are assigned using the relevant boundary conditions. That is, on perfect electric conductor (PEC) surfaces, ${\mathbf{\hat n}} \times \left[\kern-0.15em\left[ {\mathbf{E}} 
\right]\kern-0.15em\right] = 2{\mathbf{\hat n}} \times {{\mathbf{E}}^ - }$, ${\mathbf{\hat n}} \times \left[\kern-0.15em\left[ {\mathbf{H}} 
\right]\kern-0.15em\right] = 0$, and, for absorbing boundary conditions (ABCs), ${\mathbf{\hat n}} \times \left[\kern-0.15em\left[ {\mathbf{E}} 
\right]\kern-0.15em\right] = {\mathbf{\hat n}} \times {{\mathbf{E}}^ - }$, ${\mathbf{\hat n}} \times \left[\kern-0.15em\left[ {\mathbf{H}} 
\right]\kern-0.15em\right] = {\mathbf{\hat n}} \times {{\mathbf{H}}^ - }$~\cite{Hesthaven2008}. We note here that, since the performance of ABCs degrades rapidly with the angle of incidence, they are used to terminate perfectly matched layers (PMLs) wrapping around the simulation domain. Details of the implementation of PMLs can be found in~\cite{Gedney2009, Chen2019pml}.

${\mathbf{E}_k}({\mathbf{r}},t)$ and ${\mathbf{H}_k}({\mathbf{r}},t)$ are expanded using Lagrange polynomials ${\ell _i}({\mathbf{r}})$~\cite{Hesthaven2008}
\begin{align}
\label{expEt}
E_k^\nu ({\mathbf{r}},t) \simeq \sum\limits_{i = 1}^{{N_p}} {{E_k^\nu}({{\mathbf{r}}_i},t){\ell _i}({\mathbf{r}})}  = \sum\limits_{i = 1}^{{N_p}} {E_{k,i}^{\nu}(t){\ell _i}({\mathbf{r}})}\\
\label{expHt}
H_k^\nu ({\mathbf{r}},t) \simeq \mathop \sum \limits_{i = 1}^{{N_p}} {H_k^\nu}({{\mathbf{r}}_i},t){\ell _i}({\mathbf{r}}) = \mathop \sum \limits_{i = 1}^{{N_p}} H_{k,i}^{\nu}(t){\ell _i}({\mathbf{r}})
\end{align}
where $\mathbf{r}_i$, $i=1,...,N_p$, denote the locations of the interpolation nodes, $\nu \in \{ x,y,z\}$ denotes the components of vectors in the Cartesian coordinate system, and $E_{k,i}^{\nu}(t)$ and $H_{k,i}^{\nu}(t)$ are the unknown coefficients to be solved for. Substituting \eqref{expEt} and \eqref{expHt} into \eqref{strongE} and \eqref{strongH} yields~\cite{Hesthaven2008}
\begin{align}
\label{semiMaxwell0}
& {\varepsilon _k}{\partial _t}\bar E_k^x(t) \! = \! \bar D_k^y\bar H_k^z(t) \! - \! \bar D_k^z\bar H_k^y(t) \! - \! {\bar{F}_k}F_k^{{\mathbf{E}},x}(t) \! - \! \bar J_k^x(t)\\
\label{semiMaxwell1}
& {\varepsilon _k}{\partial _t}\bar E_k^y(t) \! = \! \bar D_k^z\bar H_k^x(t) \! - \! \bar D_k^x\bar H_k^z(t) \! - \! {\bar{F}_k}F_k^{{\mathbf{E}},y}(t) \! - \! \bar J_k^y(t)\\
\label{semiMaxwell2}
& {\varepsilon _k}{\partial _t}\bar E_k^z(t) \! = \! \bar D_k^x\bar H_k^y(t) \! - \! \bar D_k^y\bar H_k^x(t) \! - \! {\bar{F}_k}F_k^{{\mathbf{E}},z}(t) \! - \! \bar J_k^z(t)\\
\label{semiMaxwell3}
& {\mu _k}{\partial _t}\bar H_k^x(t) \! = \! - [\bar D_k^y\bar E_k^z(t) \! - \! \bar D_k^z\bar E_k^y(t)] \! + \! {\bar{F}_k}F_k^{{\mathbf{H}},x}(t)\\
\label{semiMaxwell4}
& {\mu _k}{\partial _t}\bar H_k^y(t) \! = \! - [\bar D_k^z\bar E_k^x(t) \! - \! \bar D_k^x\bar E_k^z(t)] \! + \! {\bar{F}_k}F_k^{{\mathbf{H}},y}(t)\\
\label{semiMaxwell5}
& {\mu _k}{\partial _t}\bar H_k^z(t) \! = \! - [\bar D_k^x\bar E_k^y(t) \! - \! \bar D_k^y\bar E_k^x(t)] \! + \! {\bar{F}_k}F_k^{{\mathbf{H}},z}(t)
\end{align}
where $\bar E_k^\nu(t)  = {[E_{k,1}^{\nu}(t),...,E_{k,{N_p}}^{\nu}(t)]^T}$ and $\bar H_k^\nu(t)  = {[H_{k,1}^{\nu}(t),...,H_{k,{N_p}}^{\nu}(t)]^T}$ are unknown vectors, $\bar J_k^\nu(t)  = {[J_{k,1}^{\nu}(t),...,J_{k,{N_p}}^{\nu}(t)]^T}$ is the current density vector, $J_{k,i}^{\nu}(t) = {J_k^\nu }({{\mathbf{r}}_i},t)$, and $F_k^{{\mathbf{E}},\nu }(t)$ and $F_k^{{\mathbf{H}},\nu }(t)$ are the $\nu$ components of 
$\hat{\mathbf{n}} \times [{\mathbf{H}_k}(t) - {\mathbf{H}}_k^*(t)]$ and $\hat{\mathbf{n}} \times [{\mathbf{E}_k}(t) - {\mathbf{E}}_k^*(t)]$,
respectively. $\bar D_k^\nu = {\bar{M}_k^{-1}} \tilde S_k^\nu $ and $\bar{F}_k = {\bar{M}_k^{-1}}  {\tilde L_k}$, where ${\bar M_k}$ is the mass matrix defined the same as before, ${\tilde S_k}$ and ${\tilde L_k}$ are the stiffness matrix and the surface mass matrix defined as
\begin{align*}
& \tilde S_k^\nu (i,j) = \int_{{\Omega _k}} {{\ell _i}({\mathbf{r}}) \partial_{\nu}{\ell _j (\mathbf{r})} } dV\\
& {\tilde L_k}(i,j) = \oint_{\partial {\Omega _k}} {{\ell _i}({\mathbf{r}}){\ell _j}({\mathbf{r}})dS}
\end{align*}
respectively. We note here that $\varepsilon ({\mathbf{r}})$ and $\mu ({\mathbf{r}})$ are assumed constant in each element. The semi-discretized system \eqref{semiMaxwell0}-\eqref{semiMaxwell5} is integrated in time using a low-storage five-stage fourth-order Runge-Kutta method~\cite{Hesthaven2008}.

\subsection{Comments}
\label{comments}
Several comments about the formulation of the multiphysics solver described in the previous sections are in order. (i) The formulation of the DG scheme developed for solving/updating the Maxwell equations ignores the frequency dependency of the dielectric permittivity for the sake of simplicity in description. But also this is a good approximation since the optical EM wave excitation has a very narrow bandwidth (less than $1\%$). We should note here that the frequency dependency of the dielectric permittivity can easily be accounted for using several well-known methods as described in~\cite{Gedney2012dispersive, Li2015resistive}. (ii) {\color{black}The generation rate is simply assumed proportional to the photon flux (the absorbed power density of the optical EM wave divided by the photon energy corresponding to the center frequency).} Although this is widely used in simulation of optoelectronic devices~\cite{Chuang2012, Khabiri2012, Moreno20141, Burford2016, Bashirpour2017}, a more strict treatment should consider the frequency dependency of the absorption of the semiconductor {\color{black} and the photon energy}. Nevertheless, since the bandwidth of the optical EM wave is narrow and it dominates the THz EM wave in magnitude, this non-dispersive simplified generation model performs very well in real devices~\cite{Chuang2012, Khabiri2012, Moreno20141, Burford2016, Bashirpour2017}. 
(iii) We should emphasize here that both of the semi-discrete systems \eqref{semiDD0}-\eqref{semiDD1} and \eqref{semiMaxwell0}-\eqref{semiMaxwell5}) are integrated in time using explicit schemes to yield the time samples of the unknown coefficients. The coupling of the solutions is carried out as explained in Section \ref{coupling}. We note here that the resulting time marching scheme does not call for inversion of any matrix systems during the time updates. Therefore, MPI-based parallelization of this time-domain solver on distributed memory systems is rather straightforward and leads to a high scaling efficiency~\cite{Chen2019discontinuous, Chen2019parallel}. 
(iv) It is clear from~\eqref{semiDD0}-\eqref{semiDD1} that $\bar{Q}_k(t)$ is a local variable which is the reason why the discretization method is termed local DG~\cite{Cockburn1998}. This variable is allocated as a temporary array of size $3{N_p}$ that is flushed repeatedly during the loop over the elements. Therefore, it does not increase the memory requirements of the solver.

\begin{table}[!t]
	\centering
	\begin{threeparttable}
		\renewcommand{\arraystretch}{1.2}
		\caption{Physical parameters used for the PCD examples}
		\label{parameters}
		\centering
		\begin{tabular}{c c}
			\hline
			\hline
			Laser\tnote{a} & $f_c=375$ THz, $f_{w}=25$ THz, Power $=0.63$ mW \\ 
			LT-GaAs & $\epsilon_r=13.26$, $\mu_r=1.0$, $\alpha=1 \mu\mathrm{m}^{-1} $, $\eta=1.0$ \\ 
			SI-GaAs & $\epsilon_r=13.26$, $\mu_r=1.0$ \\ 
			Metal\tnote{b} & $\epsilon_{\infty}=1$, $\omega_p=9.03/{\hbar}$, $\gamma=0.053q/{\hbar}$ \\ \hline
			Temperature & 300 K \\ 
			$V_{bias}$ & 10 V \\ 
			C & $1.3\times 10^{16}$ cm$^{-3}$ \\
			$n_i$ & $9\times 10^6$ cm$^{-3}$ \\ \hline
			Mobility &
			{ \begin{math} \begin{array}{cc}
				\mu_e^0=8000 \mathrm{cm}^2\mathrm{/V/s}, \mu_h^0=400 \mathrm{cm}^2\mathrm{/V/s} \\
				V_e^{sat}\!=\!1.725 \! \times \! 10^{7} \mathrm{cm/s}, V_h^{sat} \! = \! 0.9 \! \times \! 10^{7}\mathrm{cm/s} \\
				\beta_e=1.82, \beta_h=1.75
				\end{array} \end{math} } \\ \hline
			Recombination & 
			{ \begin{math} \begin{array}{cc}
				\tau_e=0.3 \mathrm{ps}, \tau_h=0.4 \mathrm{ps} \\
				n_{e1}=n_{h1}=4.5 \times 10^6 \mathrm{cm}^{-3} \\
				C_e^A=C_h^A=7\times 10^{-30} \mathrm{cm}^6\mathrm{/s}
				\end{array} \end{math} } \\
			\hline
			\hline
		\end{tabular}
		\smallskip
		\scriptsize
		\begin{tablenotes}
			\item[a] {$f_c$ is the center frequency and $f_w$ is the bandwidth}
			\item[b] {Drude model parameters~\cite{Olmon2012}}
		\end{tablenotes}
	\end{threeparttable}
\end{table}

\begin{figure}[!tbp]
	\centerline{\includegraphics[width=0.6\columnwidth]{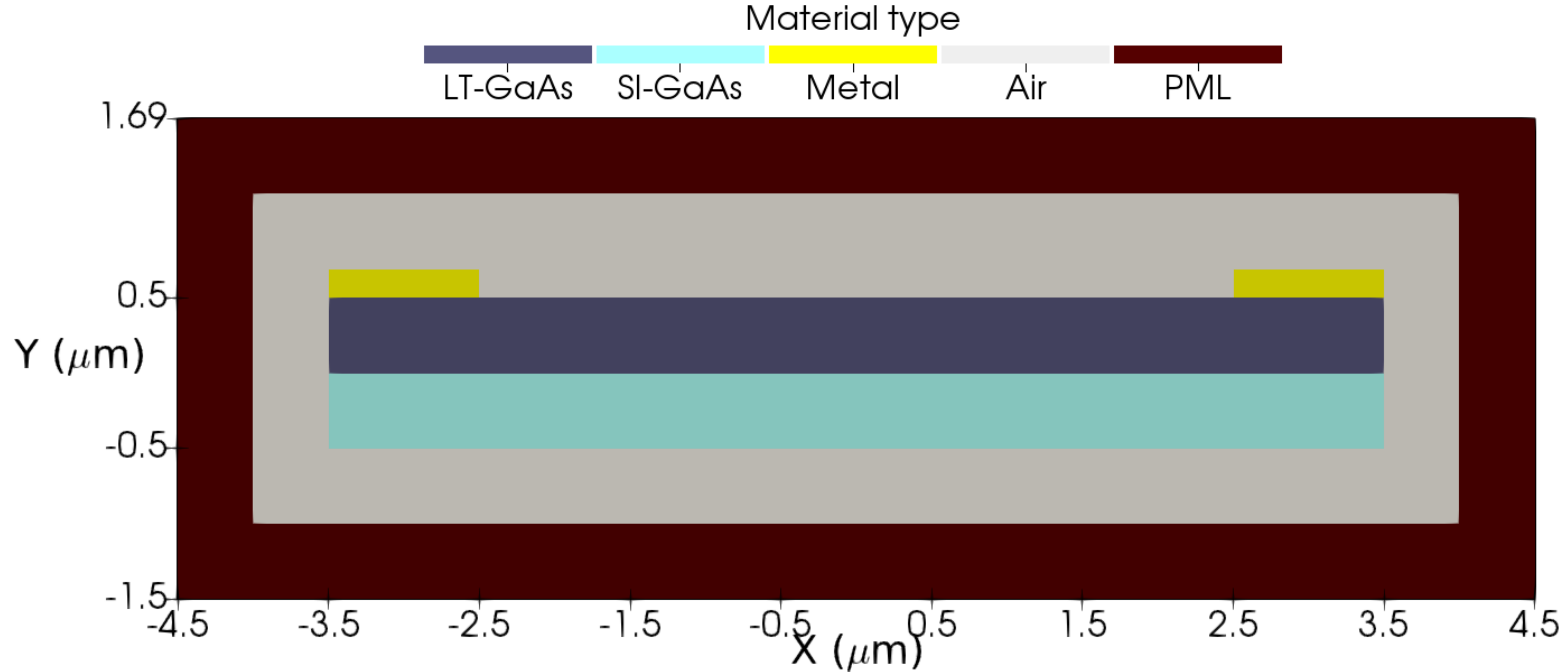}}
	\caption{Cross section of the conventional PCD considered in the first example.}
	\label{PCA2D}
\end{figure}

\begin{figure}[!bhp]
	\centering
	\subfloat[\label{PCA2D_SSa}]{\includegraphics[width=0.65\columnwidth]{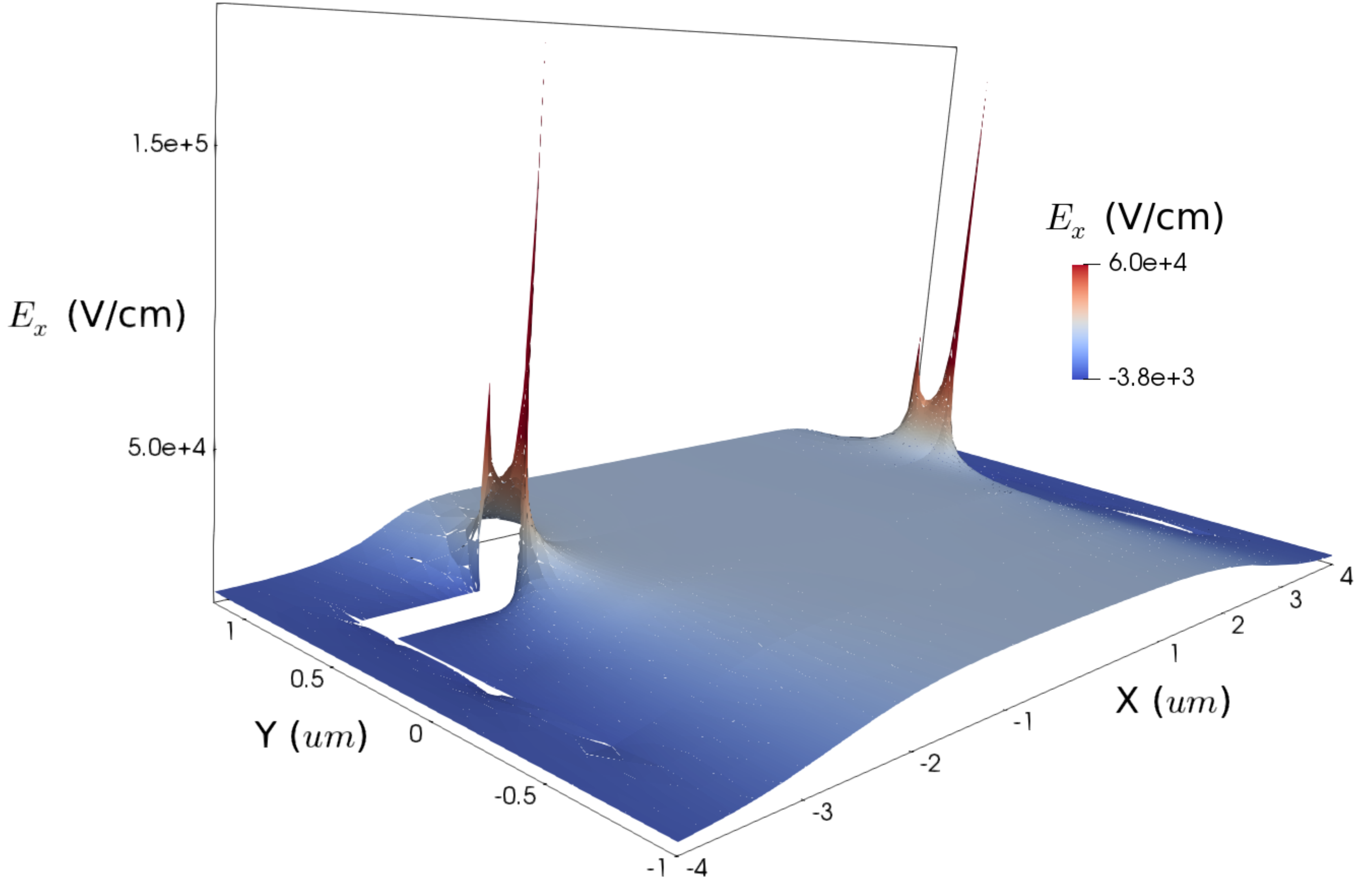}}\\
	\subfloat[\label{PCA2D_SSb}]{\includegraphics[width=0.65\columnwidth]{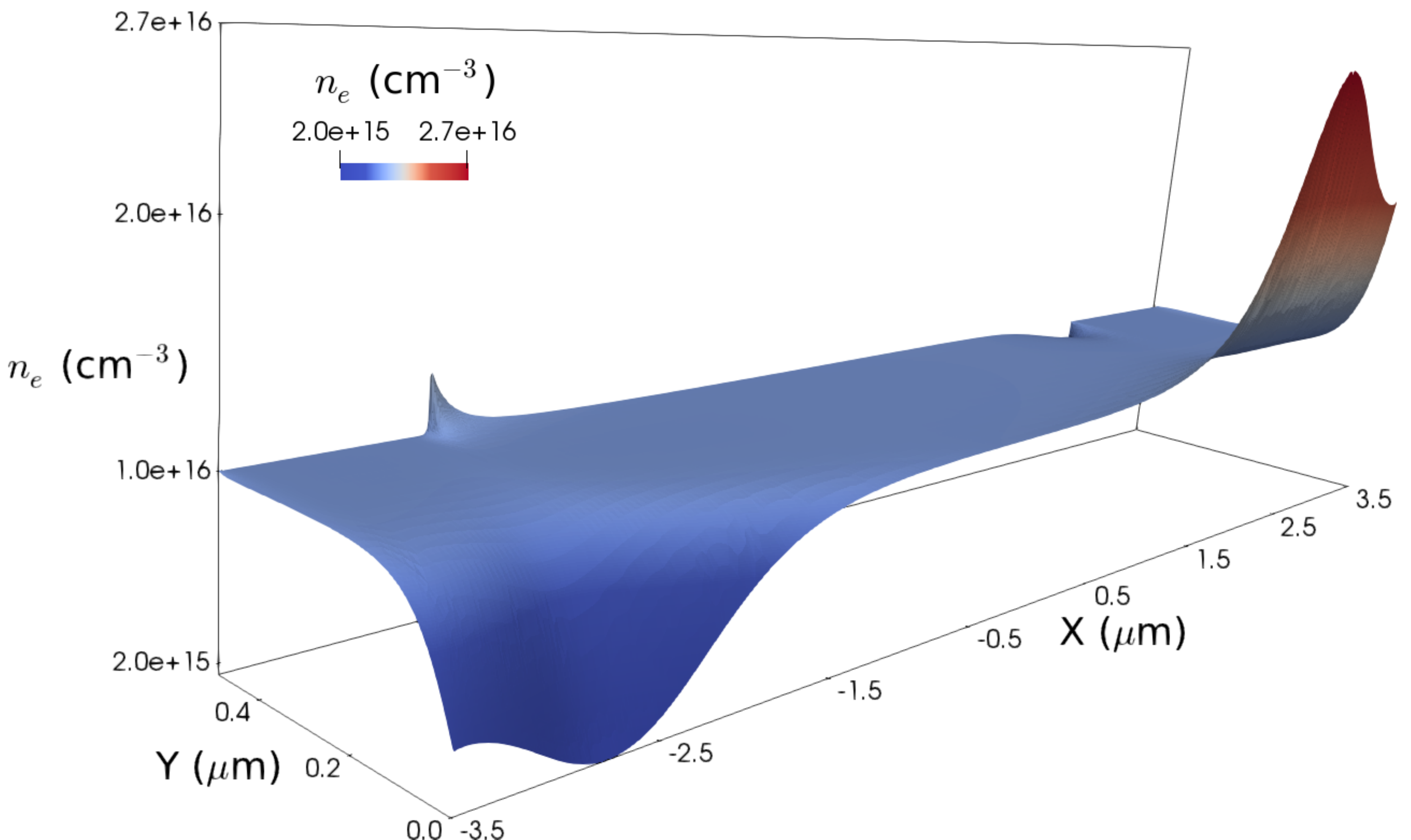}}\\
	\subfloat[\label{PCA2D_SSc}]{\includegraphics[width=0.65\columnwidth]{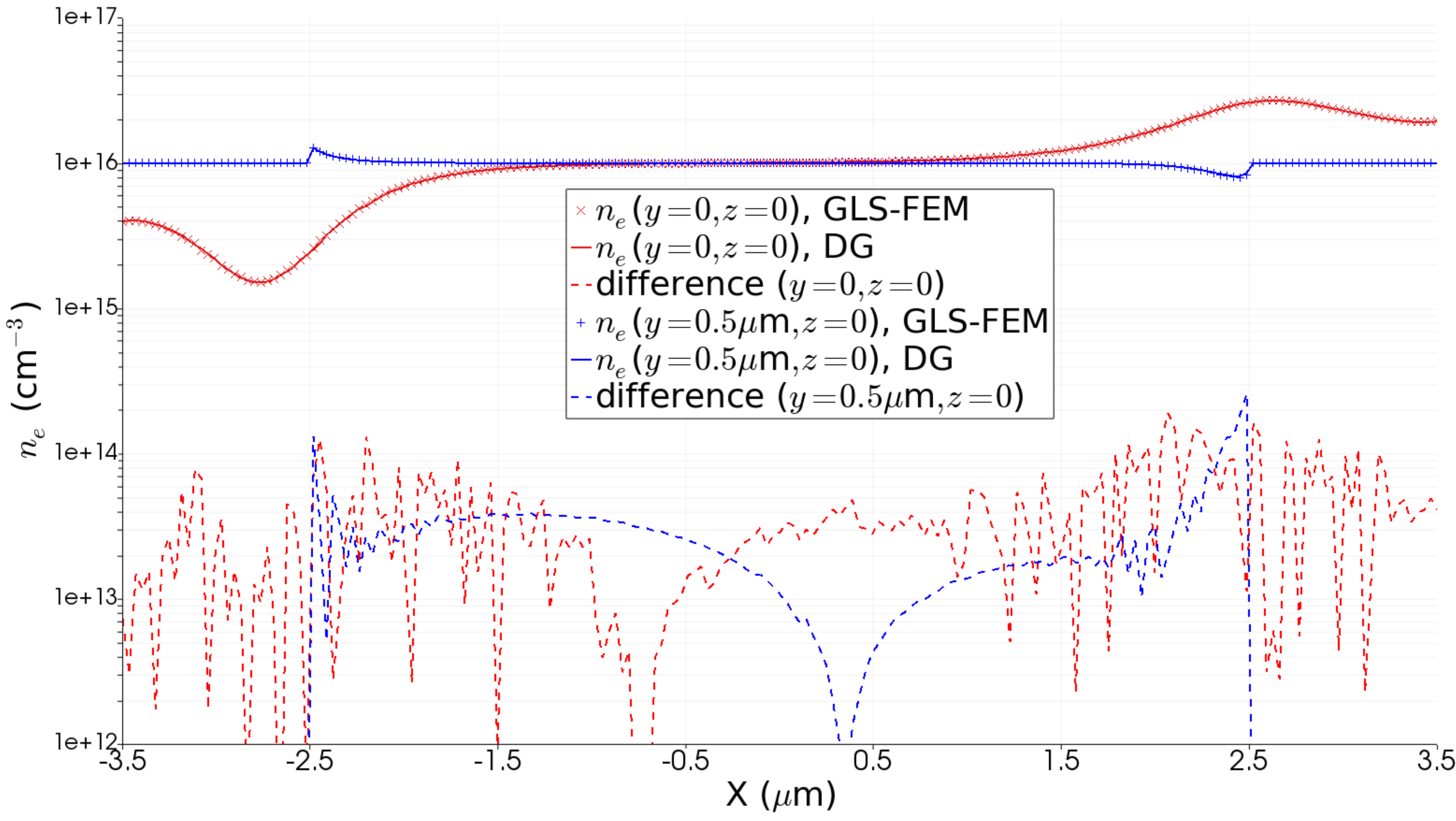}}
	\caption{(a) Stationary electric field and (b) electron density computed using the proposed DG scheme. (c) Electron density computed using the DG-based framework and COMSOL (GLS-FEM) along the lines $(y=0, z=0)$ and $(y=0.5~\mu \mathrm{m}, z=0)$.}
	\label{PCA2D_SS}
\end{figure}

\begin{figure*}[!h]
	\centering
	\subfloat[\label{PCA2D_TSa}]{\includegraphics[width=0.75\columnwidth]{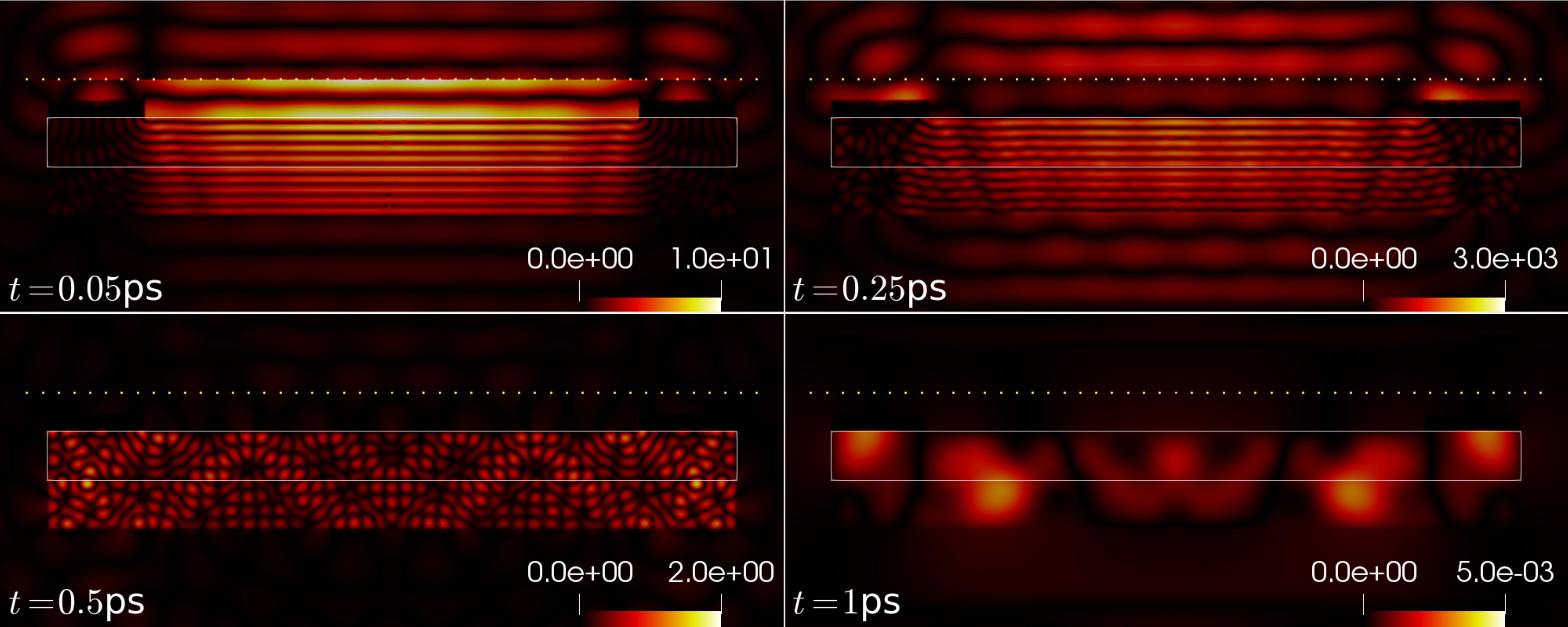}}
	\hspace{0.01cm}
	\subfloat[\label{PCA2D_TSb}]{\includegraphics[width=0.5\columnwidth]{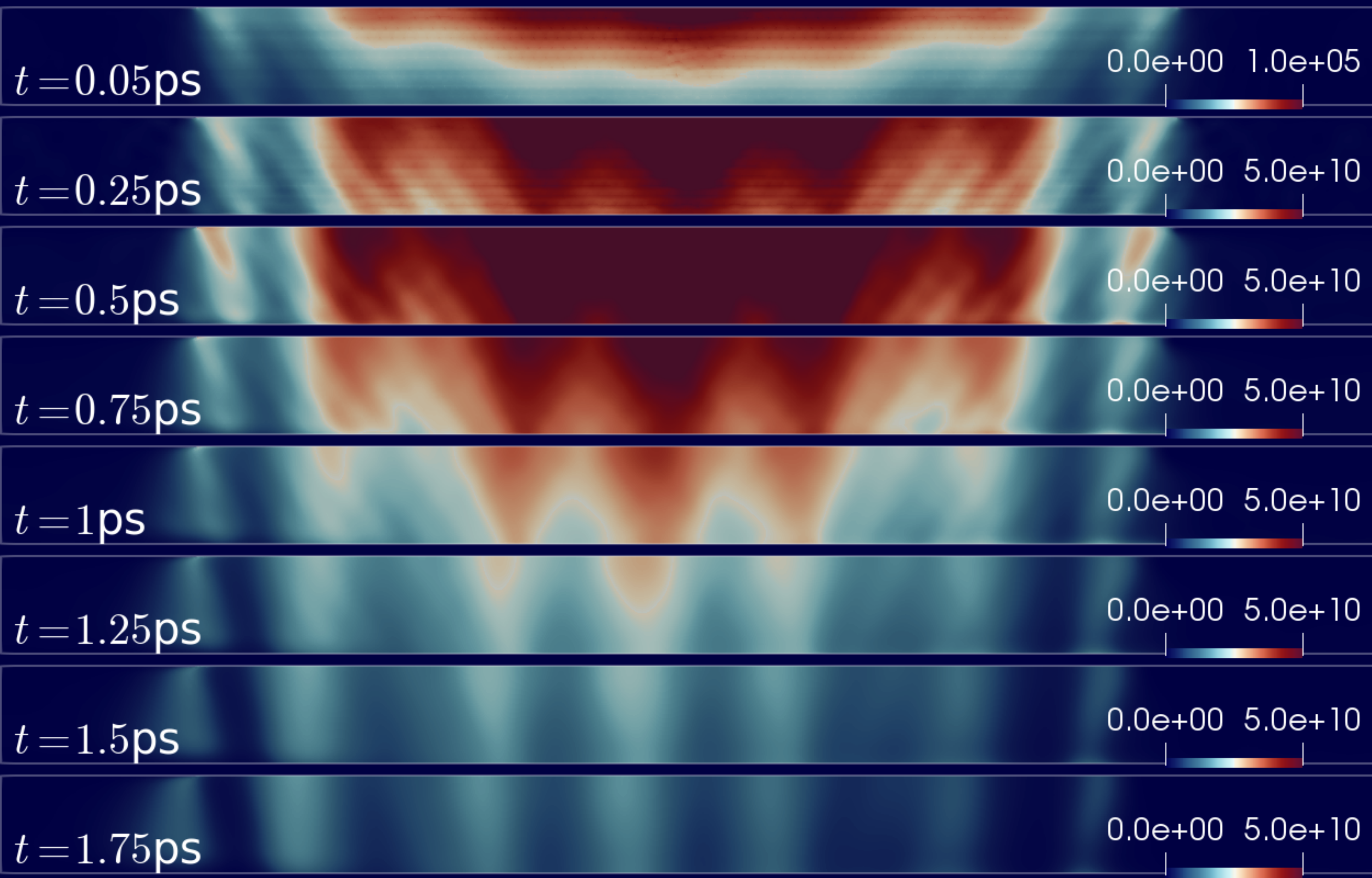}}
	\caption{(a) $|{H_z}(\mathbf{r},t)|$ computed using the DG scheme at different time instants. The yellow dotted line shows the aperture of optical EM wave excitation. The gray box indicates the LT-GaAs region. (b) $n_e(\mathbf{r},t)$ computed using the DG scheme at different time instants.}
	\label{PCA2D_TS}
\end{figure*}
\section{Numerical Results}
\label{results}
\subsection{Conventional PCD}
To demonstrate the accuracy of the proposed framework, we first use it to simulate a conventional face-to-face dipole PCD~\cite{Tani1997}. This device is selected for validation since its structure is relatively simple, and the FDTD and FEM-based approaches, which model the carrier generation due to the optical EM wave excitation using closed-form analytical expressions~\cite{Moreno2014, Burford2016, Bashirpour2017}, produce results that match well with experimental ones.

As is done in~\cite{Burford2016, Bashirpour2017}, we characterize the current response and therefore only consider the central gap region of the device. The cross section of the device is shown in Fig.~\ref{PCA2D}.  The width and height of the semiconductor and substrate layers are $7~\mu$m and $0.5~\mu$m, respectively. Both electrodes {\color{black} are made of gold and} have a width of $1~\mu$m and a height of $0.19~\mu$m. {\color{black} For the transient simulation, the permittivity of gold is represented using the Drude model~\cite{Gedney2012}. The Drude model parameters of gold at the optical frequencies~\cite{Olmon2012} as well as other physical parameters are shown in Table~\ref{parameters}.}

The physical parameters are listed in Table~\ref{parameters}. The stationary and time-dependent DD equations are solved only within the photoconductive semiconductor (LT-GaAs) layer while the Poisson equation and the time-dependent Maxwell equations are solved in the whole domain. The following boundary conditions are used.

(i) Stationary and time-dependent DD equations: On electrode/semiconductor interfaces, metal contacts are assumed ideal Ohmic contacts and Dirichlet boundary condition ${n_e} = (C + \sqrt {{C^2} + 4{n_i}^2} )/2$, ${n_h} = {n_i}^2/{n_e}$, which is derived using the local charge neutrality~\cite{Schroeder1994, Vasileska2010}, is used. On semiconductor/insulator interfaces, no carrier spills happen, which corresponds to the homogeneous Robin boundary condition~\cite{Schroeder1994, Chen2020efficient} $\hat{ \mathbf{n} } \cdot { \mathbf{J}_{e(h)}} = 0$.

(ii) Poisson equation: On electrode surfaces, assuming ideal Ohmic contact, Dirichlet boundary condition $\varphi  = {V_{external}} + {V_T} \ln{(n_e^s/n_i)}$ is enforced~\cite{Vasileska2010}. On truncation boundary of the computation domain, homogeneous Neumann boundary condition $\hat{\mathbf{n}} \cdot \nabla \varphi  = 0$  is enforced under the assumption that the static fields reaching the boundary are small and do not change the solution near the device.

(iii) Maxwell equations: The computation domain is ``wrapped around” by PMLs~\cite{Gedney2009} truncated by first-order ABCs~\cite{Hesthaven2008}.

Several comments about the space and time characteristic scales are in order.

(i) Space scales: Characteristic space scales, more specifically, the EM wavelength at the optical and THz frequencies, {\color{black}the skin depth of the EM wave in the metal $l_d$~\cite{Olmon2012}}, Debye length $l_D$~\cite{Vasileska2010}, and the Peclet number~\cite{Trangenstein2013} constrain the edge length used in the mesh discretizing the computation domain. The Debye length $l_D = \sqrt{2\varepsilon V_T \left/ \left(qn_{c} \right) \right. }$ is a measure of the carrier density variation in space~\cite{Vasileska2010}. Here, $ n_{c}$ is an estimate of maximum achievable carrier density. The Peclet number $\Delta_d C_P$, where $\Delta_d$ is the largest edge length of the local mesh element and $C_P = |\mu_c \mathbf{E}|/(2D_c)$, has to be less than $1$ to ensure the stability of the convection-diffusion component in the DD equations~\cite{Trangenstein2013}. We should note here that the potential distribution is smooth, therefore the mesh using the edge length determined by the space scales described above discretizes it accurately. {\color{black} The values of these space-scale parameters for the materials used in this example (see Table~\ref{parameters}) and the corresponding required edge lengths are provided in Table~\ref{length}.} We should note that the last number in Table~\ref{length} is the smallest geometry dimension and the mesh used for this example represents the geometry very accurately. 

(ii) Time scales: The Courant-Friedrichs-Lewy (CFL) stability conditions of Maxwell equations~\cite{Hesthaven2008} and diffusion and drift components of the DD equations~\cite{Liu2016, Wang2015},  constrain the time-step sizes used in the transient solution. Table~\ref{time} lists the values of the maximum time-step sizes allowed by these three CFL conditions with the edge length obtained from the characteristic space scale discussion above. We should note that the time-step sizes required to resolve the periods of the optical and THz EM waves are larger than those required by their CFL conditions. Table~\ref{time} clearly shows that the DD equations can be integrated using a larger time-step size than the one required for the integration of the Maxwell equations. 
\begin{table}[!t]
	\renewcommand{\arraystretch}{1.15}
	\caption{Length scales}
	\label{length}
	\centering
	\begin{tabular}{c c c}
		\hline
		\hline
		Quantity & Value (nm) & Required mesh size (nm)\\
		\hline
		optical wavelength & $800$ & $\sim 200$\\
		THz wavelength & $3\times10^5$ & $\sim 10^5$\\
		$l_d$ & $25$ & $\sim 10$\\
		$l_D$ & $30$ & $\sim 10$\\
		$C_P^{-1}$ & $10$ & $\sim10$\\
		geometrical details & $100$ & $\sim 100$\\
		\hline
		\hline
	\end{tabular}
\end{table}

\begin{table}[!t]
	\renewcommand{\arraystretch}{1.15}
	\caption{Time scales}
	\label{time}
	\centering
	\begin{tabular}{c c}
		\hline
		\hline
		Quantity & Required time-step size (ps)\\
		\hline
		CFL (Maxwell equations) & $\sim 10^{-7}$\\
		CFL (DD drift term) & $\sim 10^{-6}$\\
		CFL (DD diffusion term) & $\sim 10^{-6}$\\
		\hline
		\hline
	\end{tabular}
\end{table}

Using the discussion above on the characteristic space scales as a guide, the semiconductor layer {\color{black}and the metallic electrodes} are discretized using a mesh with a minimum edge length of  $10~\mathrm{nm}$ and the maximum edge length of the mesh is allowed to reach $70~\mathrm{nm}$ in the GaAs layers and $200~\mathrm{nm}$ in the rest of the computation domain. Similarly, following the discussion on the characteristic time scales, the time-step sizes for the Maxwell and DD equations are selected to be ${10^{-7}}~\mathrm{ps}$ and $5 \times {10^{-7}}~\mathrm{ps}$, respectively.

\begin{figure}[!t]
	\centering
	\subfloat[\label{PCA2D_Currenta}]{\includegraphics[width=0.6\columnwidth]{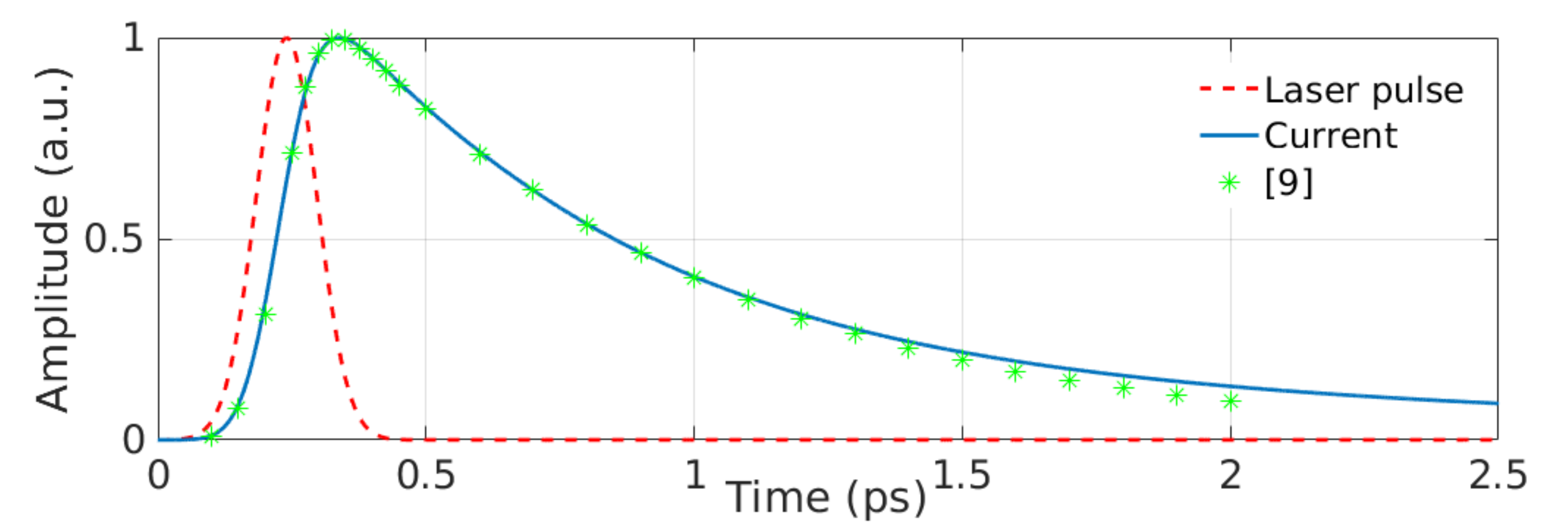}}\\
	\subfloat[\label{PCA2D_Currentb}]{\includegraphics[width=0.6\columnwidth]{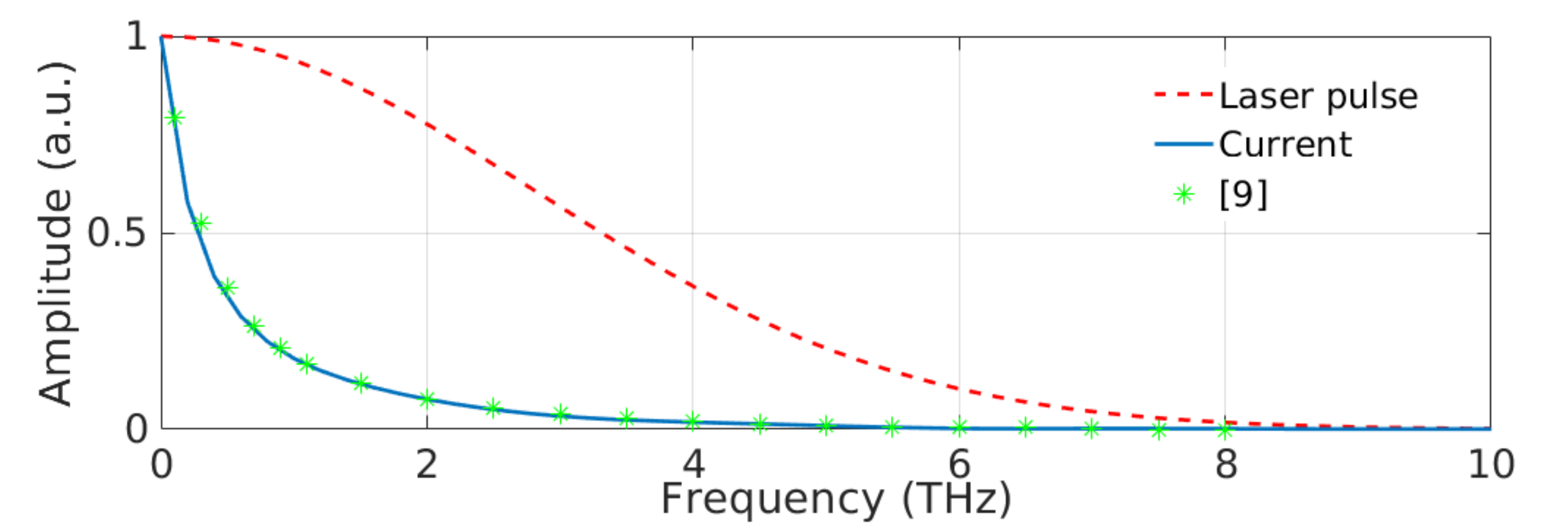}}
	\caption{(a) The time signature and (b) the spectrum of the optical EM wave excitation and the generated THz current.}
	\label{PCA2D_Current}
\end{figure}

Figs.~\ref{PCA2D_SS} (a) and (b) show the stationary state electric field and electron density computed using the DG-based framework on the plane $z=0$, respectively. Fig.~\ref{PCA2D_SS} (b) shows that the electron density distribution has sharp boundary layers, which is a typical phenomenon observed in the solution of convection dominated convection-diffusion problems. A finer mesh is needed near the boundaries as the Peclet number gets larger. Fig.~\ref{PCA2D_SS} (c) compares the electron density computed using the DG scheme and the COMSOL semiconductor module~\cite{COMSOL} along the lines $(y=0, z=0)$ and $(y=0.5~\mu m, z=0)$. Note that for the COMSOL module, Galerkin least-square (GLS)-FEM is used for stabilization and we choose the fully-coupled nonlinear scheme and the second order discretization~\cite{COMSOL}. The “difference” shown in Fig.~\ref{PCA2D_SS} (c) is computed using $\left\| n_e^{\mathrm{DG}}-n_e^{\mathrm{FEM}} \right\|$, where $ n_e^{\mathrm{DG}}$ and $n_e^{\mathrm{FEM}}$ refer to the solutions computed by the DG scheme and the COMSOL. The figure shows that the difference is two orders of magnitude smaller than the solutions demonstrating the accuracy of stationary state DG solution.

\begin{figure}[!t]
	\centerline{\includegraphics[width=0.6\columnwidth]{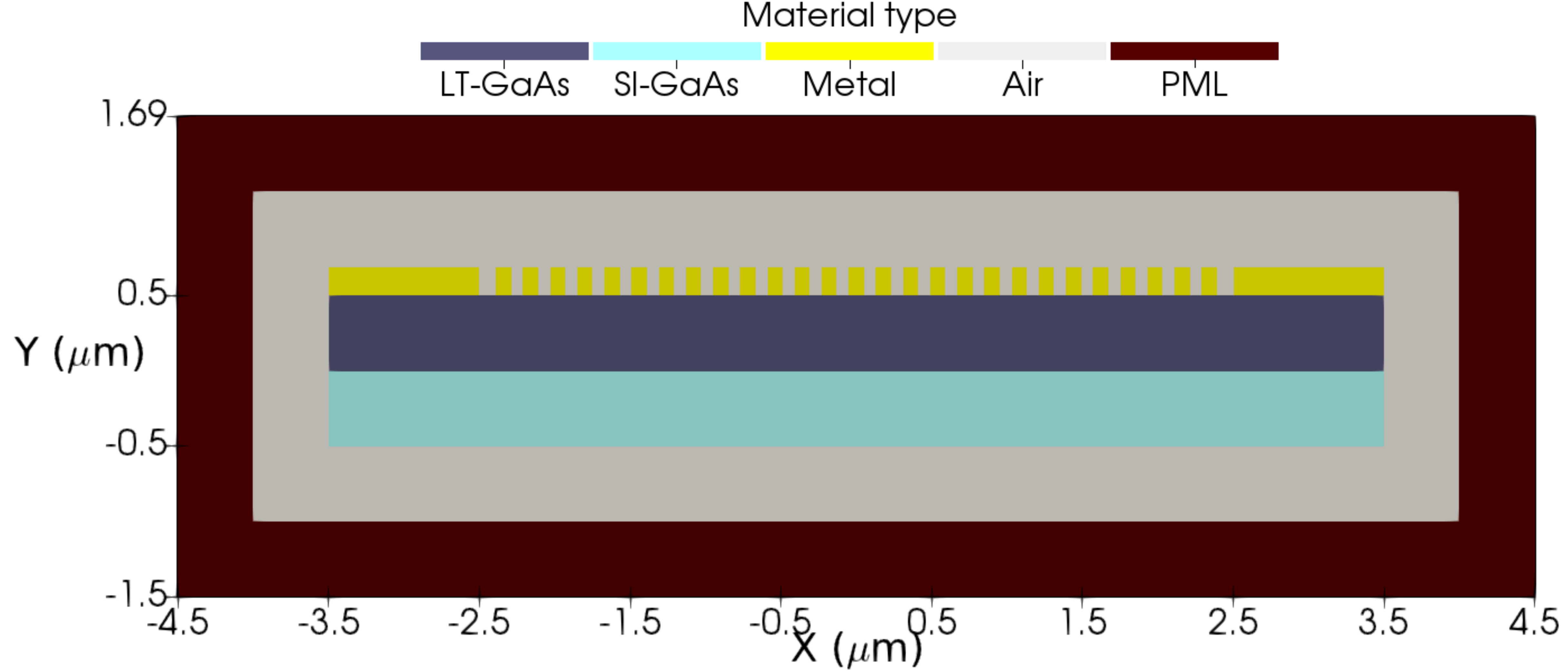}}
	\caption{Cross section of the plasmonic PCD.}
	\label{PlasmPCA2D}
\end{figure}

\begin{figure}[!t]
	\centering
	\subfloat[\label{PlasmPCA2D_SSa}]{\includegraphics[width=0.65\columnwidth]{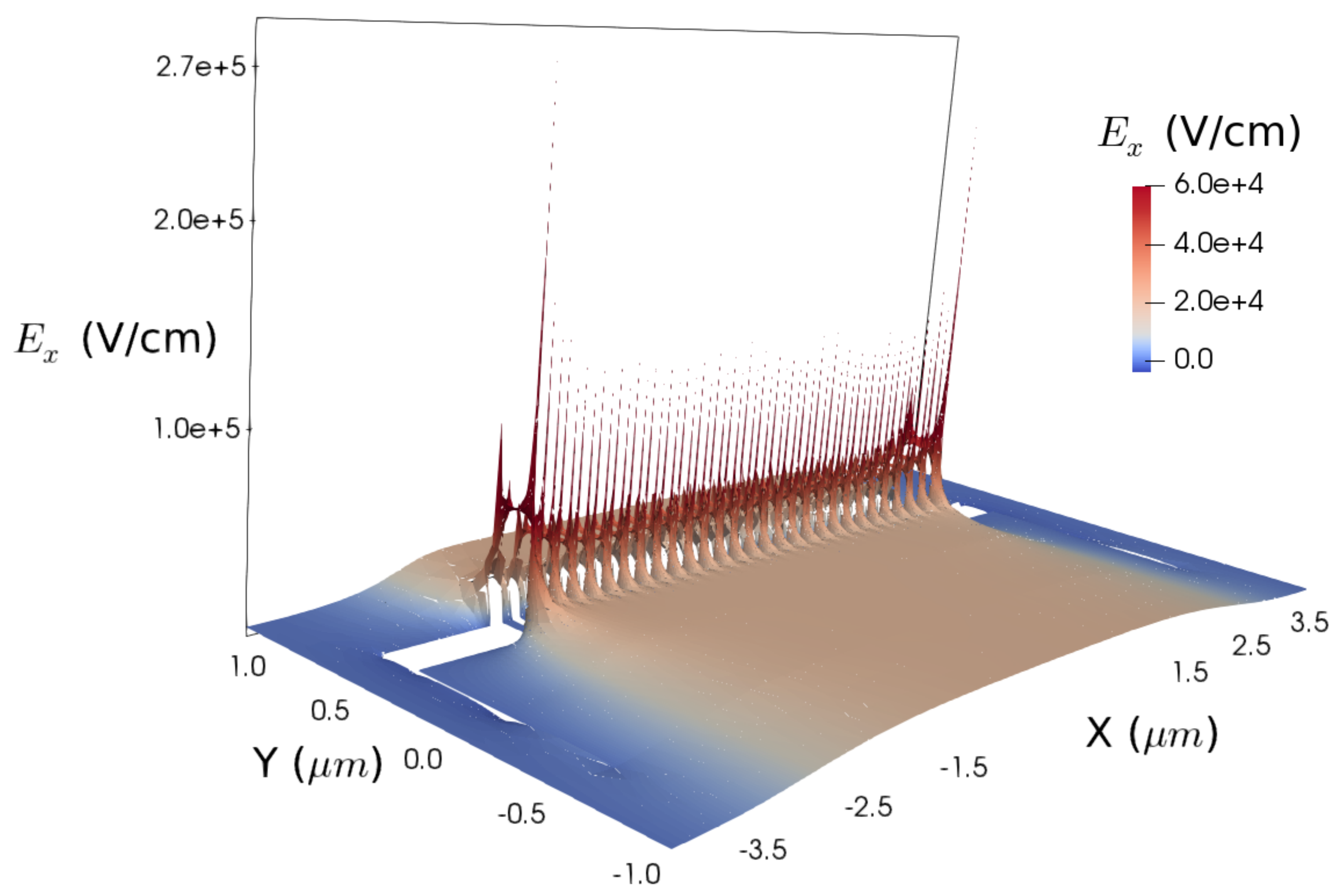}}\\
	\subfloat[\label{PlasmPCA2D_SSb}]{\includegraphics[width=0.65\columnwidth]{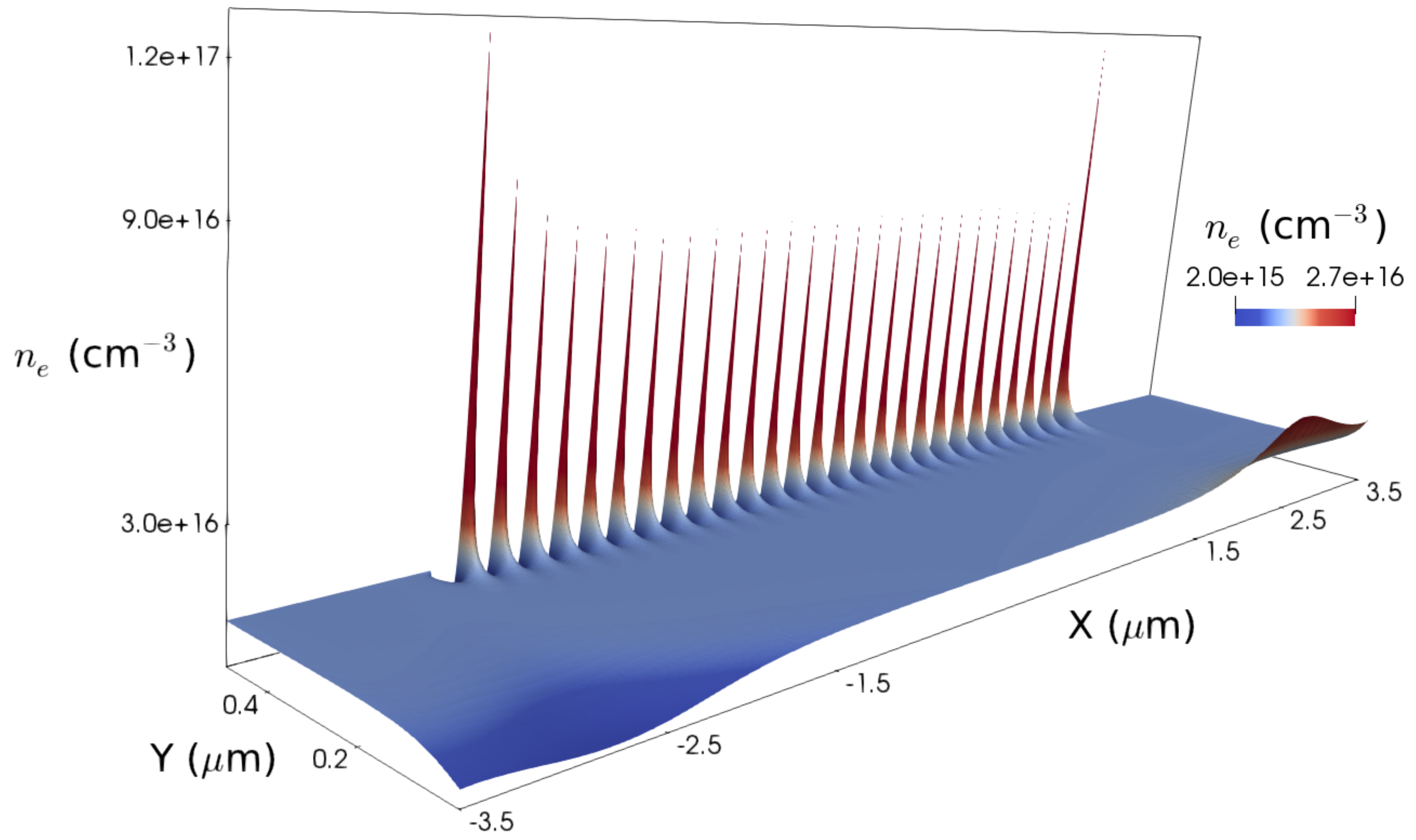}}\\
	\subfloat[\label{PlasmPCA2D_SSc}]{\includegraphics[width=0.65\columnwidth]{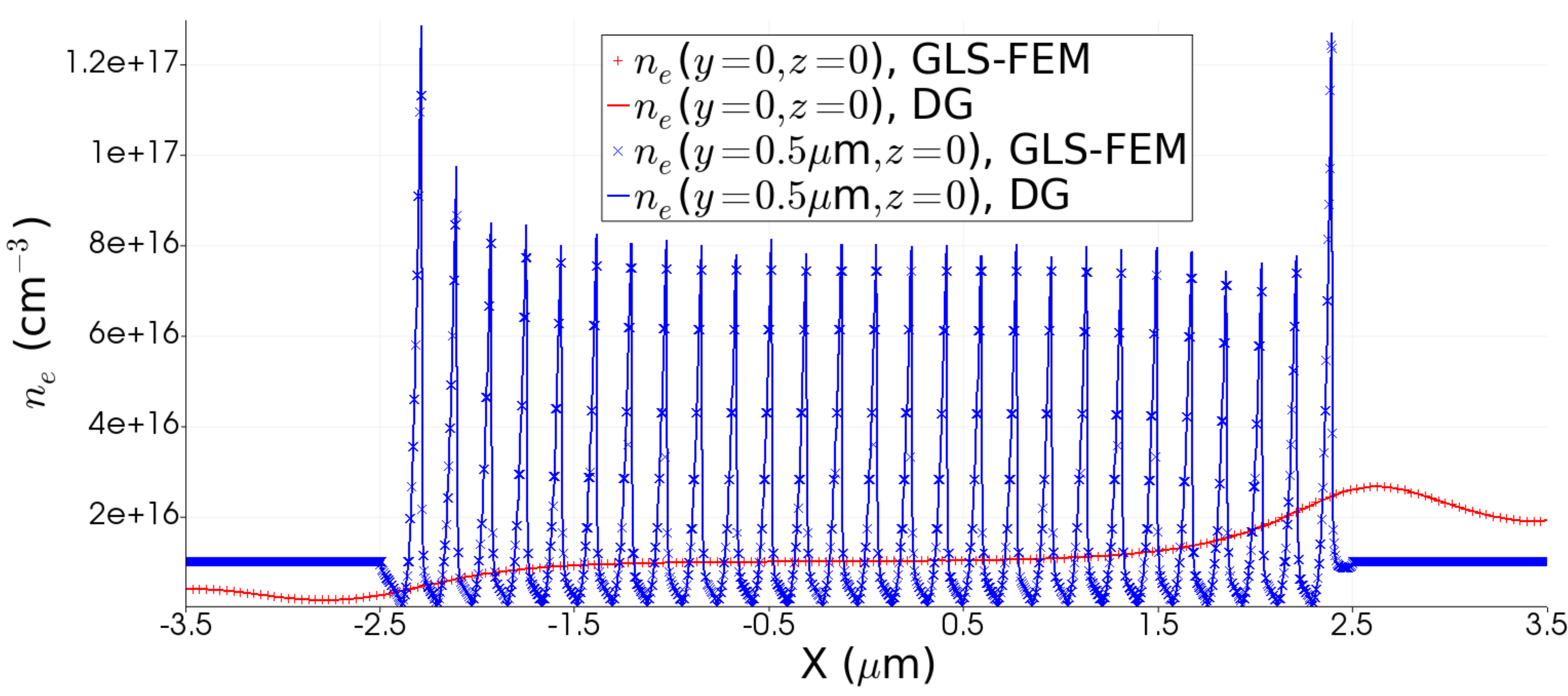}}
	\caption{(a) Stationary electric field and (b) electron density computed using the proposed DG scheme. (c) Electron density computed using the DG-based framework and COMSOL (GLS-FEM) along the lines $(y=0, z=0)$ and $(y=0.5~\mu m, z=0)$.}
	\label{PlasmPCA2D_SS}
\end{figure}

\begin{figure*}[!t]
	\centering
	\subfloat[\label{PlasmPCA2D_TSa}]{\includegraphics[width=0.75\columnwidth]{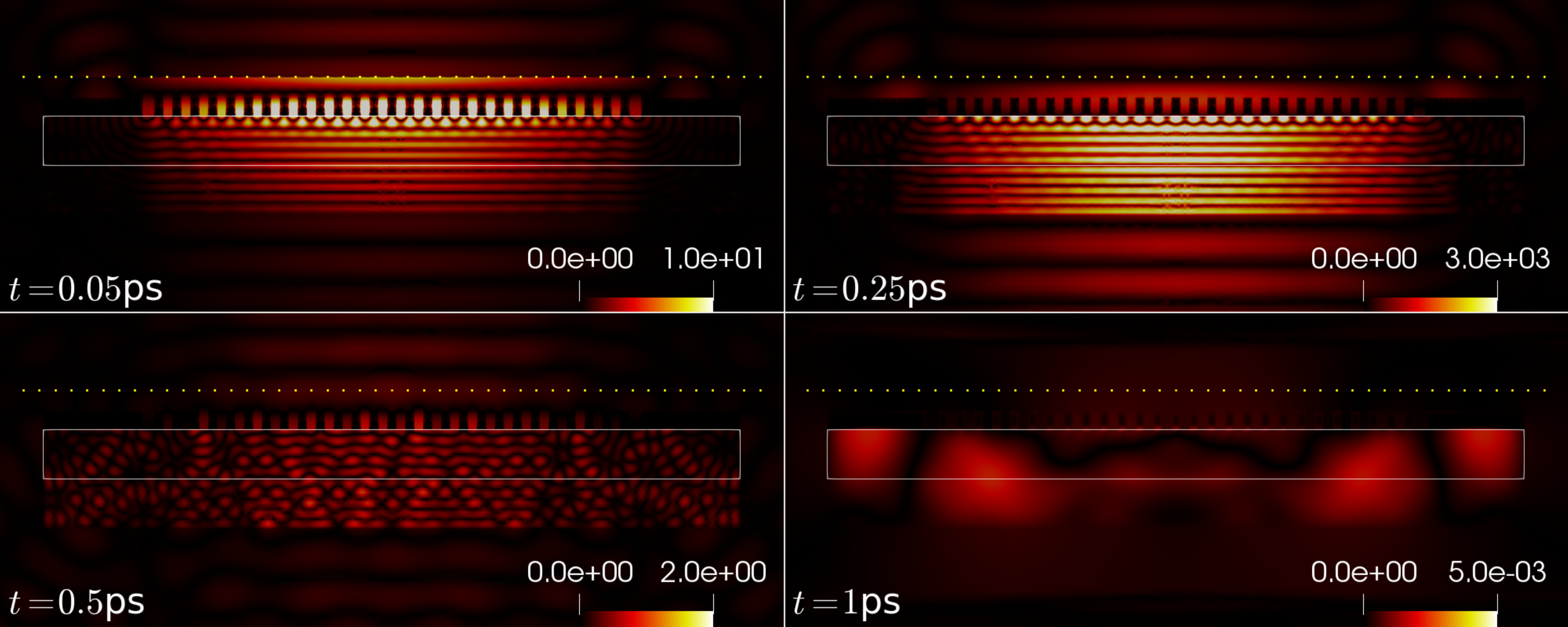}}
	\hspace{0.01cm}
	\subfloat[\label{PlasmPCA2D_TSb}]{\includegraphics[width=0.5\columnwidth]{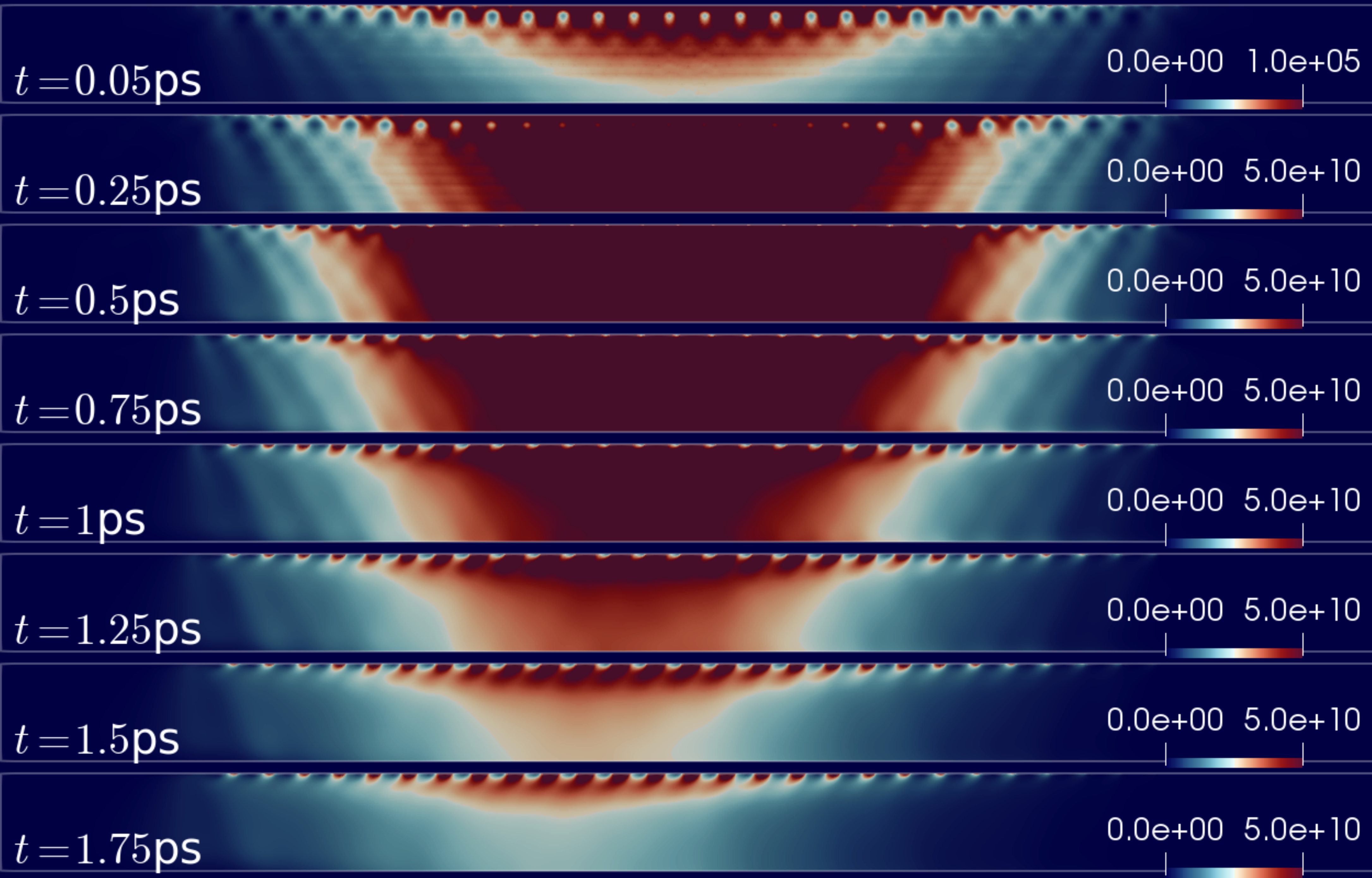}}\\
	\caption{(a) $|{H_z}(\mathbf{r},t)|$ computed using the DG scheme at different time instants. The yellow dotted line shows the aperture of optical EM wave excitation. The gray box indicates the LT-GaAs region. (b) $n_e(\mathbf{r},t)$ computed using the DG scheme at different time instants.}
	\label{PlasmPCA2D_TS}
\end{figure*}
Once the stationary state response of the device is obtained, multiple transient simulations can be executed with different excitations. The aperture of the optical EM wave source is located at $y=0.8~\mu\mathrm{m}$, generating optical EM waves propagating from top (aperture) to bottom (device). The pulse shape parameters are shown in Table~\ref{parameters}. The intensity of the EM field on the aperture has a Gaussian distribution with beam width $3~\mu\mathrm{m}$.

Figs.~\ref{PCA2D_TS} (a) and (b) show $\left| H_z^t(\mathbf{r},t) \right|$ and $n_e^t(\mathbf{r},t)$ at several time instants, respectively, and Figs.~\ref{PCA2D_Current} (a) and (b) show the time signatures of the optical EM wave excitation and the generated THz current and their spectrum obtained using Fourier transform. Figs.~\ref{PCA2D_TS} (a) and (b) and Fig.~\ref{PCA2D_Current} (a) also serve as a reference for the following discussion.

At $t=0.05~\mathrm{ps}$, the optical EM wave generated by the aperture arrives the semiconductor layer. Because  the permittivity of LT-GaAs is high, a large part of the incident field’s energy is reflected back (see the reflected wave above the air-semiconductor interface and behind the aperture). At the same time, in the LT-GaAs layer, the incident field’s energy entering the device is partially absorbed and carriers are generated near the air-semiconductor interface. At $t=0.25~\mathrm{ps}$, the incident field reaches its pulse peak and the electron density increases to $\sim10^{11}~\mathrm{cm}^{-3}$. The short incident field pulse passes quickly and, after $t=0.5~\mathrm{ps}$, only some scattered fields reside in the high permittivity region due to internal reflections. During that time, the electron density keeps increasing until the excitation pulse decays to $20$\% of its peak value (at $t\approx0.4~\mathrm{ps}$, see Fig.~\ref{PCA2D_Currenta}. After $t\approx0.4~\mathrm{ps}$, $n_e^t(\mathbf{r},t)$ decays slowly due to the recombination.

Comparing the electron density distributions at different time instants, it can be clearly observed that electrons move toward the anode on the left side. The picture of holes (not shown) is similar but with holes moving toward the cathode on the right side with a lower speed. The resulting current shown in Figs.~\ref{PCA2D_Current} (a) and (b) match very well the result presented in~\cite{Moreno2014}. Because of the simplicity of the single interface scattering, the recorded EM field intensity in the semiconductor layer has almost the same pulse shape as the optical excitation signal. This explains why the analytical generation rate~\cite{Moreno2014} works very well.

\subsection{Plasmonic PCD}

Next, to demonstrate the applicability of the proposed DG-based framework, we use it to simulate a plasmonic PCD. The cross section of the device is shown in Fig.~\ref{PlasmPCA2D}. For this example, gold nanostructures are added between the two electrodes. The periodicity, the duty cycle, and the thickness of the nanograting are $180~\mathrm{nm}$, $5/9$, and $190~\mathrm{nm}$, respectively.

The surfaces of the nanostructures are modeled as floating potential surfaces for the stationary state simulation (i.e., for solution of the Poisson equation)~\cite{Chen2020steadystate, Chen2020float}. Note that in the case of nanostructures being used as electrodes~\cite{Yang2014}, one could use Dirichlet boundary condition with the bias voltage on the surface of the nanostructure. We should also note here that a finer mesh is used near and inside the nanostructures to correctly resolve the geometry and the exponential decay of the plasmonic fields. {\color{black} In this region, the mesh has a minimum edge length of $3~\mathrm{nm}$ and the maximum edge length of the mesh is allowed to reach $70~\mathrm{nm}$ in the GaAs layers and $200~\mathrm{nm}$ in the rest of the computation domain.} The corresponding time-step size used for the time integration of the Maxwell equations is $0.3\times10^{-7}~\mathrm{ps}$. Other simulation parameters and boundary conditions are same as those in the previous example.

Figs.~\ref{PlasmPCA2D_SS} (a) and (b) show the stationary state electric field and electron density computed using the DG-based framework on the plane $z=0$, respectively. Significant enhancement of the static electric field and carrier densities is observed near the grating-semiconductor interfaces. Because of velocity saturation (see~\cite{Vasileska2010,Chen2020steadystate}), the enhanced static field results in a significant drop of the mobility, which further influences the transient response. Fig.~\ref{PlasmPCA2D_SS} (c) compares the electron density computed using the DG scheme and the COMSOL semiconductor module~\cite{COMSOL} along the lines $(y=0, z=0)$ and $(y=0.5~\mu m, z=0)$. Excellent agreement is observed between the results.

\begin{figure}[!t]
	\centering
	\subfloat[\label{PlasmPCA2D_Currenta}]{\includegraphics[width=0.55\columnwidth]{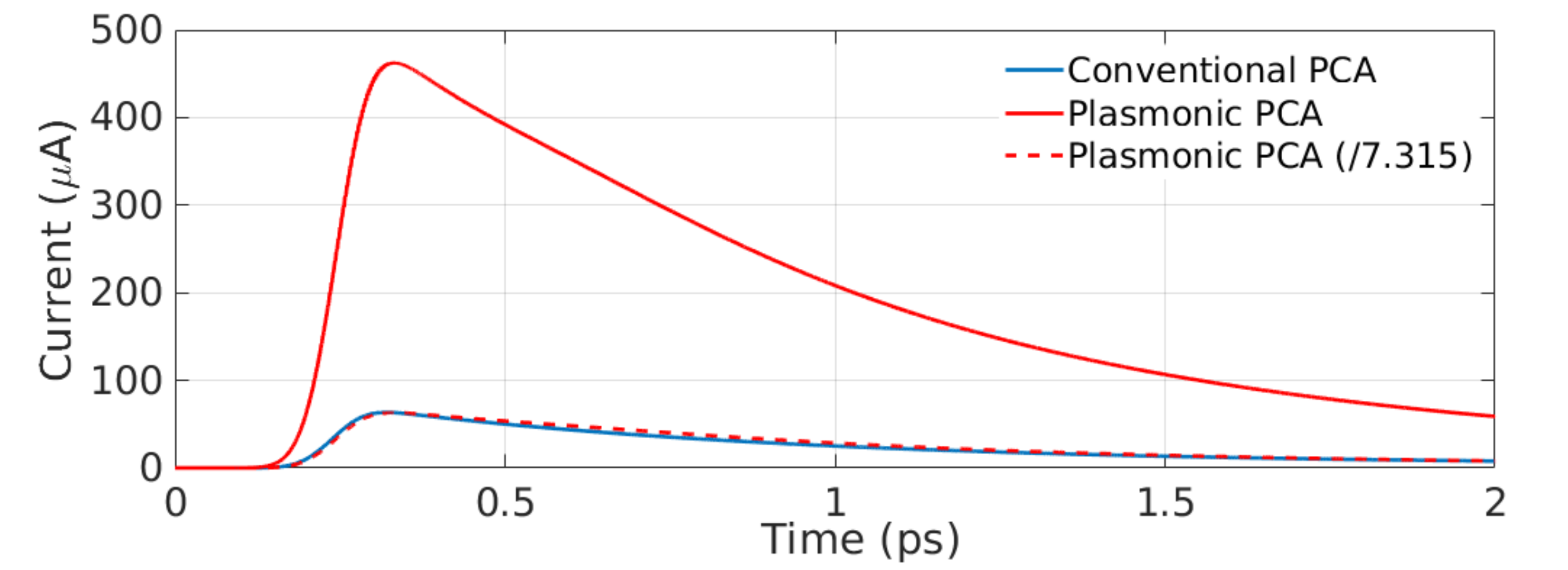}}\\
	\subfloat[\label{PlasmPCA2D_Currentb}]{\includegraphics[width=0.55\columnwidth]{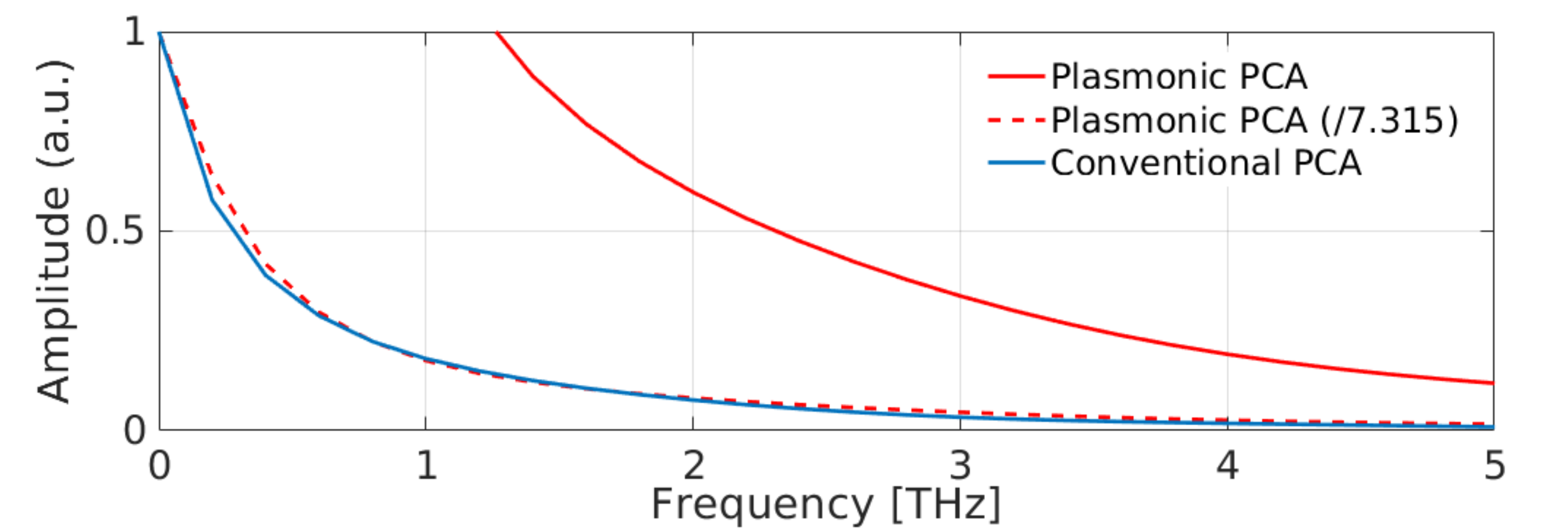}}
	\caption{(a) The time signature and (b) the spectrum of the THz current generated on the plasmonic PCD and on the conventional PCD.}
	\label{PlasmPCA2D_Current}
\end{figure}
Figs.~\ref{PlasmPCA2D_TS} (a) and (b) show $\left| H_z^t(\mathbf{r},t) \right|$ and $n_e^t(\mathbf{r},t)$ at several time instants, respectively. For comparison, Fig.~\ref{PlasmPCA2D_TS} uses the same color scale as Fig.~\ref{PCA2D_TS} for the previous example.

A comparison of $\left| H_z^t(\mathbf{r},t) \right|$ shown in Fig.~\ref{PlasmPCA2D_TS} and Fig.~\ref{PCA2D_TS} shows that the transient EM fields on the plasmonic PCD are much stronger. In addition, plasmonic mode patterns are observed at $t=0.05~\mathrm{ps}$ and $t=0.25~\mathrm{ps}$. Accordingly, the level of electron density in Fig.~\ref{PlasmPCA2D_TS} is much higher and shows an inhomogeneous pattern. As time marching goes on, electrons drift toward the anode. Figs.~\ref{PlasmPCA2D_Current} (a) and (b) compare the time signatures of the current induced on the plasmonic PCD and the conventional one (previous example) and their spectrum obtained using Fourier transform. The figures clearly show that the inclusion of the nanostructures in the PCD design enhances the current by almost $7$ times .

\section{Conclusion}
\label{conclusion}
A multiphysics framework is developed for the simulation of plasmonic PCDs. The device is modeled in two stages: (i) A stationary stage that corresponds to the nonequilibrium state under a bias voltage and is mathematically described by the nonlinearly coupled system of Poisson and DD equations and (ii) a transient stage that corresponds to the dynamic EM field and carrier interactions under an optical EM wave excitation and is mathematically described by the nonlinearly coupled system of the time-dependent Maxwell and DD equations, respectively. This paper focuses on the latter. More specifically, we develop a DG scheme to discretize the time-dependent Maxwell and DD equations. The resulting semi-discrete systems are integrated in time using explicit schemes (with different time-step sizes) to yield the samples of the unknowns (EM field intensities and carrier densities). During the time integration, the nonlinear coupling between the two sets is accounted for by alternately feeding updated solutions into each other. The discretization flexibility provided by the DG and the fact that two different time-step sizes are used to integrate the two systems, respectively, helps to accurately and efficiently account for the widely different characteristic scales of the Maxwell and DD equations. 

The accuracy of the multiphysics framework is verified against that of previously developed methods via simulations of a conventional PCD. A plasmonic PCD is simulated using the proposed framework to demonstrate its capability to account for the complex physics involved in the optical-to-THz process when plasmonic metallic nanostructures are present.

The operation of plasmonic PCDs involves different mechanisms that are responsible for the overall device performance~\cite{Lepeshov2017review, Burford2017review, Kang2018review, Yardimci2018review, Yachmenev2019review, Piao2000, Loata2007}, such as the nanostructure-tailored bias electric field, the plasmonic modes, the carrier-screening effect, the high-power saturation, and the antenna radiation efficiency, etc. Studying the dependence of the device performance on these different mechanisms is a challenging task and is often expensive to do with experiments. The proposed numerical framework provides an efficient way to optimize the device performance by making it feasible to analyze different mechanisms separately by controlling/interfering the physical parameters in numerical experiments. Finally, we want to note here that extensions of the proposed framework can be developed for simulation of similar semiconductor optoelectronic devices, such as photovoltaic-effect based devices.

\section*{Acknowledgment}
This publication is based upon work supported by the King Abdullah University of Science and Technology (KAUST) Office of Sponsored Research (OSR) under Award No 2016-CRG5-2953. The authors would like to thank the King Abdullah University of Science and Technology Supercomputing Laboratory (KSL) for providing the required computational resources.



Generated by IEEEtran.bst, version: 1.14 (2015/08/26)

\end{document}